\documentclass[a4j,11pt]{amsart}

\usepackage[mathscr]{eucal}
\usepackage{amscd}
\usepackage{amsfonts}
\usepackage{amsmath, amsthm, amssymb}
\usepackage{latexsym}
\usepackage[dvips]{graphics}
\usepackage{amssymb} 

\theoremstyle{plain} %
\newtheorem{theorem}{\indent\sc Theorem}[section] %
\newtheorem{lemma}[theorem]{\indent\sc Lemma}
\newtheorem{corollary}[theorem]{\indent\sc Corollary}
\newtheorem{proposition}[theorem]{\indent\sc Proposition}

\theoremstyle{definition} %
\newtheorem{definition}[theorem]{\indent\sc Definition}

\newtheorem{remark}[theorem]{\indent\sc Remark}
\newtheorem{example}[theorem]{\indent\sc Example}

\def\C{{\mathbb{C}}}
\def\R{{\mathbb{R}}}
\def\Q{{\mathbb{Q}}}
\def\Z{{\mathbb{Z}}}
\def\N{{\mathbb{N}}}
\def\G{{\mathbb{G}}}%
\def\D{{\mathbb{D}}}

\title[The Gauss map of minimal surfaces]{Value distribution theoretical properties of the Gauss map of pseudo-algebraic minimal surfaces}

\author[Y.~Kawakami]{Yu Kawakami}

\address{
Graduate school of Mathematics, 
Nagoya University, 
Nagoya, 464-8602
/Japan
}
\email{m02008w@math.nagoya-u.ac.jp}

\begin{document}
\maketitle
\begin{abstract}
In this thesis, we study value distribution theoretical properties of the Gauss map of 
pseudo-algebraic minimal surfaces in $n$-dimensional Euclidean space. 
After reviewing basic facts, we give estimates for the number of exceptional values and 
the totally ramified value numbers and the corresponding unicity theorems for them. 
\end{abstract}

\tableofcontents

\section{Introduction}
In this thesis we study some value distribution theoretical properties of the Gauss map of complete minimal 
surfeces.

In $3$-dimensional Euclidean space $\R^{3}$, the Gauss map of a minimal surface is considered as a meromorphic 
function on the corresponding open Riemann surface. This has given rise to many problems, because via 
Enneper-Weierstrass representation there is remarkable analogy between the minimal surface theory and the value 
distribution theory of meromorphic functions. Among them, a very interesting problem is the following : 
if $g$ is the Gauss map of a non-flat complete regular minimal surfece, how many values can $g$ omit? 
In 1961 Osserman \cite{O2} proved that $g$ can omit at most a subset of points of logarithmic capacity $0$. 
In 1981, Xavier \cite{X} proved that $g$ can omit at most $6$ values. Finally, in 1988 Fujimoto \cite{F1} proved that 
$g$ can omit at most $4$ values. Since there are a lot of examples of complete minimal surfaces whose Gauss maps omit $4$ values 
(for example, the classical Scherk surface), ``4'' is the best possible upper bound. Moreover, Fujimoto \cite{F4} proved that the totally ramified value number 
$\nu_{g}$, which gives more detailed information than the number of exceptional values $D_{g}$, satisfies $\nu_{g}\leq 4$, and 
this inequality is also best possible. Since the totally ramified value number is a rational number, 
it is remarkable that its upper bound is an integer.

On the other hand, Osserman \cite{O1} proved that the Gauss map of non-flat algebraic minimal surfaces $M$ can omit 
at most $3$ values. By an algebraic minimal surface, we mean a complete regular minimal surface with finite total curvature. 
However there is no known example with $D_{g}=3$. On the other hand there are many examples, of almost all topological types, 
with $D_{g}=2$ (\cite{SF} and \cite{MS}). Therefore, it has been widely believed that the sharp upper bound of $D_{g}$ is ``$2$". 
Moreover, as in the case of Fujimoto's theorem, it has been implicitly believed that the same is true for the totally 
ramified value number $\nu_{g}$. 

In this situation, we obtained the following results. At first, we investigate the totally ramified value number of 
algebraic minimal surfaces and discovered algebraic minimal surfaces with $\nu_{g}=2.5$, i.e., strictly larger than $2$ \cite{Ka}. 
This overthrew the above implicitly believed upper bound ``$2$''. Next, the author, Kobayashi and Miyaoka \cite{KKM} introduced 
the class of pseudo-algebraic minimal surfaces in $\R^{3}$ which includes algebraic minimal surfaces as a proper subclass and obtained estimates for the 
number of exceptional values and the totally ramified value number of the Gauss map in this class. These are the best possible 
estimates for pseudo-algebraic minimal surfaces and some special cases of algebraic minimal surfaces. By these estimates, 
we can understand the relationship between Fujimoto's result and Osserman's in this class and reveal the geometrical 
meaning behind them. We also give a kind of unicity theorem, which asks the least number of values at which if two Gauss maps 
$g_{1}$ and $g_{2}$ have the same inverse image then $g_{1}=g_{2}$, for this class. 
Moreover we \cite{Ka2} also proved an analog of these results for pseudo-algebraic minimal surfaces in $4$-dimensional Euclidean space $\R^{4}$, 
revealed the relationship between Fujimoto's result \cite{F1} and Hoffman-Osserman's \cite{HO1} for this class. Recently, Jin and Ru 
\cite{JR} extend our results to algebraic minimal surfaces in $n$-dimensional Euclidean space $\R^{n}$ and get a ramification estimate for this class. 

In this thesis, we give a detailed review of all the above mentioned, i.e., we study value distribution theoretical properties of the Gauss map of 
pseudo-algebraic minimal surfaces in Euclidean space. We believe our results give an important step toward the solution of Osserman's question 
(i.e., $D_{g}\leq 2$? for algebraic minimal surfaces) and a new perspective to study the Gauss map of minimal surfaces.

This paper is organized as follow. In Section 2, we recall basic facts on complete minimal surfaces in Euclidean space 
used in this paper. In Section 3, we first explain the Gauss map of minimal surfaces in $\R^{3}$ and define the pseudo-algebraic 
minimal surfaces in $\R^{3}$. Next, we give some results of value distribution theoretical properties (ramification estimate 
and unicity theorem) for this class. Finally, we give some results on Nevanlinna theory of the Gauss map. 
In Section 4, we give some analogous results (ramification estimate and unicity theorem) for pseudo-algebraic 
minimal surfaces in $\R^{4}$. 
In Section 5, we extend Jin and Ru's result and give some results on the Gauss map of pseudo-algebraic minimal surfaces in $\R^{n}$.

\bigskip
\textbf{Acknowledgement.}  \\
The author thanks Profs Ryoichi Kobayashi, Shin Nayatani and Reiko Miyaoka for 
many helpful comments and suggestions.
\section{Basic facts on complete minimal surfaces in $\R^{n}$}
In this section, we collect some basic facts on complete minimal surfaces, in particular algebraic minimal 
surfaces. 
\subsection{The Gauss map of minimal surfaces in $\R^{n}$} 
At first, we recall some facts on a surface in Euclidean space.  Let $M$ be a oriented real 2-dimensional differentiable manifold and $x=(x^{1},\ldots, x^{n})\colon M \to \R^{n}$ 
is an immersion. For a point $p\in M$, take a local coordinate system $(u^{1}, u^{2})$ around $p$ which are positively oriented. The tangent 
plane of $M$ at $p$ is given by 
\[
T_{p}(M)=\{\lambda \frac{\partial x}{\partial u^{1}} + \mu \frac{\partial x}{\partial u^{2}}\: |\: \lambda, \mu\in\R\}
\]
and the normal space of $M$ at $p$ is given by 
\[
N_{p}(M)=\{N\: |\Bigl(N, \frac{\partial x}{\partial u^{1}}\Bigr)=\Bigl(N, \frac{\partial x}{\partial u^{2}}\Bigr)=0\}
\]
where $(X, Y)$ denotes the inner product of vectors $X$ and $Y$. The metric $ds^{2}$ on $M$ induced from the standard 
metric on $\R^{n}$, called the first fundamental form on $M$, is given by 
\[
ds^{2} = (dx, dx) 
       = g_{11}(du^{1})^{2}+2g_{12}du^{1}du^{2}+g_{22}(du^{2})^{2}
\]
where
\[
g_{ij}=\Bigl(\frac{\partial x}{\partial u^{i}}, \frac{\partial x}{\partial u^{j}}\Bigr) \quad (1\leq i, j \leq 2)
\]
and the second fundamental form of $M$ with respect to a unit normal vector $N$ is given by 
\[
d\sigma^{2} 
       = b_{11}(N)(du^{1})^{2}+2b_{12}(N)du^{1}du^{2}+b_{22}(N)(du^{2})^{2}
\]
where
\[
b_{ij}(N)=\Bigl(\frac{\partial^{2} x}{\partial u^{i}\partial u^{j}}, N\Bigr) \quad (1\leq i, j \leq 2)\:.
\]
Then the mean curvature of $M$ for the normal direction $N$ at $p$ is defined by 
\begin{equation}\label{meancurvature} 
H_{p}(N)=\frac{g_{11}b_{22}(N)+g_{22}b_{11}(N)-2g_{12}b_{12}(N)}{2\{g_{11}g_{22}-(g_{12})^{2}\}}\;.
\end{equation}
\begin{definition}\label{minimalsurface}
A surface $M$ is called a {\em minimal surface} in $\R^{n}$ if $H_{p}(N)=0$ for all $p\in M$ and $N\in N_{p}(M)$.
\end{definition}
A local coordinate system $(u^{1}, u^{2})$ on an open set $U$ in $M$ is called {\em isothermal} on $U$ if 
$ds^{2}$ can be represented as 
\[
ds^{2}=\lambda^{2}\{(du^{1})^{2}+(du^{2})^{2}\}
\]
for a positive smooth function $\lambda$ on $U$. 
\begin{theorem}[S.~S.~Chern, \cite{Cr}]\label{Chern}
For every surface $M$, there is a system of isothermal local coordinates whose domains cover the whole $M$.
\end{theorem}

\begin{proposition}\label{iso}
For an oriented surface $M$ with a metric $ds^{2}$, if we take two positively oriented isothermal local coordinates 
$(u, v)$ and $(x, y)$, then $w=u+\sqrt{-1}v$ is a holomorphic function in $z=x+\sqrt{-1}y$ on the common domain.
\end{proposition}

Let $x\colon M \to \R^{n}$ be an oriented surface with a Riemannian metric $ds^{2}$. To each positive isothermal 
local coordinate system $(u, v)$ we associate the complex function $z=u+\sqrt{-1}v$. By Proposition \ref{iso}, 
we may regard $M$ as a Riemann surface. Then the metric $ds^{2}$ is given by 
\[
ds^{2}={\lambda_{z}}^{2}(du^{2}+dv^{2})
\]
where
\[
{\lambda_{z}}^{2}=\Bigl(\frac{\partial x}{\partial u}, \frac{\partial x}{\partial u}\Bigr)=\Bigl(\frac{\partial x}{\partial v}, \frac{\partial x}{\partial v}\Bigr)
\]
Set complex differentiations
\[
\frac{\partial x^{i}}{\partial z}=\frac{\partial x^{i}}{\partial u}-\sqrt{-1}\frac{\partial x^{i}}{\partial v}, \;
\frac{\partial x^{i}}{\partial \bar{z}}=\Bigl(\overline{\frac{\partial x^{i}}{\partial z}}\Bigr)\;.
\]
following Osserman \cite{O2} and $\phi_{i}=(\partial x^{i}/\partial z)dz\:(i=1,\ldots,n)$, we may rewrite the metric
\begin{equation}\label{n-metric}
ds^{2}=\frac{1}{2}\Bigl(|\phi_{1}|^{2}+\ldots+|\phi_{n}|^{2}\Bigr)\:.
\end{equation}
Define the Laplacian $\triangle_{z}=\partial^{2}/\partial u^{2}+\partial^{2}/\partial v^{2}$ in terms of the complex local
coordinate $z=u+\sqrt{-1}v$. If we take another complex local coordinate $\zeta$, then we have $\triangle_{\zeta}=|dz/d\zeta|^{2}\triangle_{z}$. 
Since $\lambda_{\zeta}=\lambda_{z}|dz/d\zeta|$, the operator $\triangle=(1/{\lambda_{z}}^{2})\triangle_{z}$ does not depend on 
the choice of complex local coordinate $z$, which is called the Laplace-Bertrami operator. 
\begin{proposition}\label{Laplace}
It holds that 
\begin{enumerate}
\item[(i)] $(\triangle x, X)=0$ for each $X\in T_{p}(M)$, 
\item[(ii)] $(\triangle x, N)=2H(N)$ for each $N\in N_{p}(M)$.
\end{enumerate}
\end{proposition}
\begin{proof}
By the assumptions, we have 
\[
\lambda^{2}=\Bigl(\frac{\partial x}{\partial u}, \frac{\partial x}{\partial u}\Bigr)=\Bigl(\frac{\partial x}{\partial v}, \frac{\partial x}{\partial v}\Bigr), \: 
\Bigl(\frac{\partial x}{\partial u}, \frac{\partial x}{\partial v}\Bigr)=0
\]
Differentiating these identities, we have
\[
\Bigl(\frac{\partial^{2} x}{\partial u^{2}}, \frac{\partial x}{\partial u}\Bigr)=\Bigl(\frac{\partial^{2} x}{\partial u\partial v}, \frac{\partial x}{\partial v}\Bigr), \: 
\Bigl(\frac{\partial^{2} x}{\partial v\partial u}, \frac{\partial x}{\partial v}\Bigr)+\Bigl(\frac{\partial x}{\partial u}, \frac{\partial^{2} x}{\partial v^{2}}\Bigr)=0\:.
\]
These imply 
\[
\Bigl(\triangle_{z}x, \frac{\partial x}{\partial u}\Bigr)=\Bigl(\triangle_{z}x, \frac{\partial x}{\partial v}\Bigr)=0\:.
\]
Since $\partial x/\partial u$ and $\partial x/\partial v$ generate the tangent plane, we conclude the assertion (i) of Proposition \ref{Laplace}. 
On the other hand, for every normal vector $N$ to $M$ it holds that 
\[
H(N)=\frac{b_{11}(N)+b_{22}(N)}{2\lambda^{2}}=\frac{(\triangle x, N)}{2}\:.
\]
It shows (ii) of Proposition \ref{Laplace}.
\end{proof}
\begin{theorem}\label{minimal1}
Let $x=(x^{1},\ldots, x^{n})\colon M\to \R^{n}$ be a surface. Then $M$ is minimal if and only if 
each $x^{i}$ is a harmonic function on $M$, namely, 
\[
\triangle_{z}x^{i}=\Bigl(\frac{\partial^{2}}{\partial u^{2}}+\frac{\partial^{2}}{\partial v^{2}}\Bigr)x^{i}=0 \quad (1\leq i \leq n)
\]
for every complex local coordinate $z=u+\sqrt{-1}v$.
\end{theorem}
\begin{proof}
By (i) of Proposition \ref{Laplace}, $\triangle_{z}x=0$ if and only if $\triangle_{z}x$ is perpendicular to the 
normal space of $M$. This is equivalent to the condition $H(N)=0$ for each $N\in N_{p}(M)$ by (ii) of Proposition \ref{Laplace}.
\end{proof}
\begin{corollary}\label{topological}
There exists no compact minimal surface without boundary in $\R^{n}$.
\end{corollary}
\begin{proof}
For a minimal surface $x=(x^{1},\ldots, x^{n})\colon M\to \R^{n}$, if $M$ is compact, then each $x^{i}$ takes 
the maximum values at a point in $M$. By the maximum principle of harmonic functions, $x^{i}$ is a constant. 
This is imposible because $x$ is an immersion.  
\end{proof}
Next, we recall the  Gauss map of a surface immersed in $\R^{n}$. We consider the set of all oriented 
2-linear subspaces in $\R^{n}$ and denote it by $\G_{2, n}(\R)$. We identify it with the quadric $\Q^{n-2}(\C)$ in the ($n-1$)-dimensional complex 
projective space $\mathbb{P}^{n-1}(\C)$ as following. To each $P\in\G_{2, n}(\R)$, taking a positively oriented basis $\{X, Y\}$ of $P$ 
such that 
\begin{equation}\label{basis}
|X|=|Y|,\quad (X, Y)=0
\end{equation} 
we assign the point $\Phi(P)=\pi(X-\sqrt{-1}Y)$, where $\pi$ denotes the natural projection of $\C^{n}\backslash\{0\}$ onto 
$\mathbb{P}^{n-1}(\C)$, namely, the map which maps each $p=(w^{1},\ldots,w^{n})(\not= (0,\ldots,0) )$ to the equivalence class 
\[
(w^{1}:\cdots :w^{n})=\{(cw^{1},\ldots,cw^{n})|c\in \C\backslash \{0\}\}
\]  
As is easily seen, the value $\Phi(P)$ does not depend on the choice of a positive basis of $P$ satisfying (\ref{basis}) but does only on $P$. 
On the other hand, $\Phi(P)$ is contained in the quadric 
\[
\Q^{n-2}(\C)=\{(w^{1}:\cdots:w^{n})|(w^{1})^{2}+\cdots+(w^{n})^{2}=0 \}(\subset \mathbb{P}^{n-1}(\C))\:.
\]
In fact, for a positive basis $\{X, Y\}$ satisfying (\ref{basis}) we have 
\[
(X-\sqrt{-1}Y, X-\sqrt{-1}Y)=(X, X)-2\sqrt{-1}(X, Y)-(Y, Y)=0\:.
\]
Conversely, take an arbitrary point $Q\in \Q^{n-2}(\C)$. It is easily seen that there is a unique oriented $2$-plane $P$ 
such that $\Phi(P)=Q$. This shows that $\Phi$ is bijective. Thus the 
set of all oriented $2$-planes in $\R^{n}$ is identified with the quadric $\Q^{n-2}(\C)$. 

Now, consider a surface $x=(x^{1},\ldots, x^{n})\colon M \to \R^{n}$. For each point $p\in M$, the oriented tangent 
plane $T_{p}(M)$ is canonically identified via $\Phi$ with an element of $\tilde{\G}_{2, n}(\R)\cong \Q^{n-2}(\C)$ 
after the parallel translation which maps $p$ to the origin. 
\begin{definition}\label{GGM}
{\em The Gauss map } of a surface $M$ is defined as the map of $M$ into $\Q^{n-2}(\C)$ which 
maps each point $p\in M$ to $\Phi(T_{p}(M))$.
\end{definition}
For a positively oriented isothermal local coordinate $(u, v)$ the vectors
\[
X=\frac{\partial x}{\partial u},\: Y=\frac{\partial x}{\partial v}
\]
give a positive basis of $T_{p}(M)$ satisfying the condition (\ref{basis}). Therefore, the Gauss map $g$ is 
locally given by 
\begin{equation}\label{Gauss1}
g=\pi(X-\sqrt{-1}Y)=\Bigr(\frac{\partial x^{1}}{\partial z}:\frac{\partial x^{2}}{\partial z}:\cdots:\frac{\partial x^{n}}{\partial z}\Bigr)\:,
\end{equation}
where $z=u+\sqrt{-1}v$. We may write $g=(\phi_{1}:\cdots:\phi_{n})$ with globally 
defined holomorphic $1$-forms $\phi_{i}=(\partial x^{i}/\partial z)dz\; (1\leq i \leq n)$.
\begin{theorem}\label{minimal2}
A surface $x\colon M\to \R^{n}$ is minimal if and only if the Gauss map $g\colon M\to \mathbb{P}^{n-1}(\C)$ is holomorphic.
\end{theorem}
\begin{proof}
Assume that $M$ is minimal. We then have 
\[
\frac{\partial}{\partial \bar{z}}\Bigl(\frac{\partial x}{\partial z}\Bigr)=\triangle x=0
\] 
by Theorem \ref{minimal1}. This shows that $\partial x/\partial z$ satisfies Cauchy-Riemann's equation. Hence, 
the Gauss map $g$ is holomorphic. 

Conversely, assume that $g$ is holomorphic. For a complex local coordinate $z$ we set $f_{i}=\partial x^{i}/\partial z$ 
$(1\leq i\leq n)$. After a suitable change of indices, we may assume that $f_{n}$ has no zero. Since $f_{i}/f_{n}$ are 
holomorphic, we have
\[
\triangle_{z} x^{i}=\frac{\partial^{2}x^{i}}{\partial\bar{z}\partial z}=\frac{\partial}{\partial \bar{z}}\Bigl(\frac{f_{i}}{f_{n}}f_{n}\Bigr)
=\frac{\partial}{\partial \bar{z}}\Bigl(\frac{f_{i}}{f_{n}}\Bigr)
f_{n}+\frac{f_{i}}{f_{n}}\frac{\partial f_{n}}{\partial \bar{z}}=f_{i}\frac{1}{f_{n}}\frac{\partial f_{n}}{\partial \bar{z}}
\] 
for $i=1,2,\ldots,n$. Write
\[
\frac{1}{f_{n}}\frac{\partial f_{n}}{\partial \bar{z}}=h_{1}+\sqrt{-1}h_{2}
\]
with real-valued functions $h_{1}, h_{2}$ and take the real parts of both sides of the above equation to see 
\[
\triangle_{z} x=2\Bigl(\frac{\partial x}{\partial u}h_{1}+\frac{\partial x}{\partial v}h_{2} \Bigr)\in T_{p}(M).
\]
According to (i) of Propotion \ref{Laplace}, we obtain $(\triangle_{z} x, \triangle_{z} x)=0$ and so $\triangle_{z} x=0$. 
This implies that $M$ is a minimal surface by Theorem \ref{minimal1}.
\end{proof}

We say that a holomorphic $1$-form $\phi$ on a Riemann surface $M$ has no real periods if
\[
\Re \displaystyle \int_{\gamma} \phi =0
\]
for every cycle $\gamma \in H_{1}(M,\Z)$. If $\phi$ has no real period, then the quantity
\[
x(z)=\Re \displaystyle \int_{\gamma^{z}_{z_{0}}} \phi
\]
depends only on $z$ and $z_{0}$ for a piecewise smooth curve $\gamma^{z}_{z_{0}}$ in $M$ joining $z_{0}$ and $z$ 
and so $x$ is a well-defined function of $z$ on $M$, which we denoted by 
\[
x(z)=\Re \displaystyle \int^{z}_{z_{0}} \phi
\]
from here on. Related to Theorem \ref{minimal2}, we show here the following construction theorem of 
minimal surfaces. 
\begin{theorem}\label{minimal3}
Let $M$ be an open Riemann surface and $\phi_{1},\ldots,\phi_{n}$ a collection of holomorphic $1$-forms on $M$ such that 
they have no common zeros, no real periods and locally satisfy
\begin{equation}\label{isothermal}
{f_{1}}^{2}+\cdots+{f_{n}}^{2}=0
\end{equation}
for holomorphic functions $f_{i}$ with $\phi_{i}=f_{i}dz$. Set
\begin{equation}\label{immersion}
x^{i}=\Re \displaystyle \int_{z_{0}}^{z} \phi_{i}
\end{equation}
for an arbitrarily fixed point $z_{0}$ of $M$. Then, the surface $x=(x^{1},\ldots,x^{n})\colon M\to\R^{n}$ is a minimal surface 
immersed in $\R^{n}$ such that the Gauss map is the map $g=(\phi_{1}:\cdots:\phi_{n})\colon M\to\Q^{n-2}(\C)$ and the 
induced metric is given by 
\begin{equation}\label{metric}
ds^{2}=\frac{1}{2}(|\phi_{1}|^{2}+\cdots+|\phi_{n}|^{2})
\end{equation}
\end{theorem}
\begin{remark}
We call the condition that each $\phi_{i}$ has no real periods ``period condition''.
\end{remark}
\begin{proof}
By assumption, each $x^{i}$ is well-defined single-valued function on $M$. Consider the map 
$x=(x^{1},\ldots,x^{n})\colon M\to \R^{n}$. Since $\frac{\partial x^{i}}{\partial z}=f_{i}$, by (\ref{isothermal}) we have
\[
{f_{1}}^{2}+\cdots+{f_{n}}^{2}=\Bigl(\frac{\partial x}{\partial u}, \frac{\partial x}{\partial u}\Bigr)-2\sqrt{-1}
\Bigl(\frac{\partial x}{\partial u}, \frac{\partial x}{\partial v}\Bigr)-\Bigl(\frac{\partial x}{\partial v}, \frac{\partial x}{\partial v}\Bigr)
=0
\]
for $z=u+\sqrt{-1}v$. This gives that
\[
\Bigl(\frac{\partial x}{\partial u}, \frac{\partial x}{\partial u}\Bigr)=\Bigl(\frac{\partial x}{\partial v}, \frac{\partial x}{\partial v}\Bigr),\:
\Bigl(\frac{\partial x}{\partial u}, \frac{\partial x}{\partial v}\Bigr)=0\:.
\]
Moreover, the Cauchy-Schwarz inequality implies
\[
\displaystyle \sum_{i<j}|\frac{\partial(x^{i}, x^{j})}{\partial(u, v)}|^{2}=
\Bigl(\frac{\partial x}{\partial u}, \frac{\partial x}{\partial u}\Bigr)\Bigl(\frac{\partial x}{\partial v}, \frac{\partial x}{\partial v}\Bigr)
-\Bigl(\frac{\partial x}{\partial u}, \frac{\partial x}{\partial v}\Bigr)^{2}=
\frac{1}{4}(|f_{1}|^{2}+\ldots+|f_{n}|^{2})>0
\]
which mean that $x$ is an immersion. Then, the induced metric is given by 
\[
ds^{2}=\Bigl(\frac{\partial x}{\partial u}, \frac{\partial x}{\partial u}\Bigr)(du^{2}+dv^{2})=\frac{1}{2}(|f_{1}|^{2}+\ldots+|f_{n}|^{2})|dz|^{2}
\]
and $(u, v)$ gives isothermal coordinate for the induced metric $ds^{2}$. On the other hand, by (\ref{Gauss1}) the Gauss map 
$g$ of $M$ is given by $g=(f_{1}:\cdots:f_{n})$ with holomorphic functions $f_{i}$ and so holomorphic. 
According to Theorem \ref{minimal2}, the surface $M$ is a minimal surface. 
\end{proof}
Let $M$ be a Riemann surface with a metric $ds^{2}$ which is conformal, namely, represented as
\begin{equation}\label{conformalmetric}
ds^{2}={\lambda_{z}}^{2}|dz|^{2}
\end{equation}
with a positive smooth function $\lambda_{z}$ in term of a complex local coordinate $z$.
\begin{definition}
For each point $p\in M$ we define the {\em Gauss curvature} of $M$ at $p$ by
\begin{equation}\label{curvature}
K_{ds^{2}}=-\triangle\log \lambda_{z}\Bigl(=-\frac{\triangle_{z}\log\lambda_{z}}{\lambda_{z}^{2}}\Bigr)\:.
\end{equation}
\end{definition}
For a minimal surface $M$ immersed in $\R^{n}$, using (\ref{metric}), we find expression
\begin{equation}\label{curvature2}
K_{ds^{2}}=-\frac{4|\phi \wedge {\phi}'|^{2}}{|\phi|^{6}}
\end{equation}
where
\[
|\phi|^{2}=\displaystyle \sum_{k=1}^{n}|f_{k}|^{2},\quad |\phi\wedge \phi'|^{2}=\sum_{i<j}|f_{i}f'_{j}-f_{j}f'_{i}|^{2}\:.
\]
This implies that the curvature of a minimal surface is always nonpositive.

If a minimal surface is flat (i.e., the Gauss curvature vanishes everywhere ), then (\ref{curvature2}) implies that $f_{i}/f_{i_{0}}\equiv const.$$(1\leq i\leq n)$ 
for some $i_{0}$ with $f_{i_{0}}\not\equiv 0$ and, therefore, that the Gauss map $g$ is a constant map.
\begin{proposition}
For a minimal surface $M$ immersed in $\R^{n}$, $M$ is flat, or equivalently, the Gauss map of $M$ is 
a constant map if and only if it lies in a plane.
\end{proposition}
\begin{proof}
The Gauss map of a surface which lies in a plane is obviously a constant map. Conversely, we assume that 
the Gauss map $g=(g_{1}:\cdots:g_{n})$ is constant. This means that every tangent plane $T_{p}(M)$ of 
$M\;(p\in M)$ is perpendicular to $(n-2)$ particular linearly independent normal vectors $N_{1},\ldots,N_{n-2}$. 
We then have
\[
\Bigl(\frac{\partial x}{\partial u}, N_{k}\Bigr)=\Bigl(\frac{\partial x}{\partial v}, N_{k}\Bigr)\quad (1\leq k\leq n-2)
\]
as function in local coordinates $u$ and $v$. Therefore, each $(x, N_{k})$ is a constant function for $k=1,\ldots,n-2$ and so 
$M$ lies in a plane. 
\end{proof}

\subsection{Algebraic minimal surfaces and its Gauss map}
Next, we introduce the total curvature of a minimal surface and recall some properties of algebraic minimal surfaces. 
Let $x=(x^{1},\ldots,x^{n}):M\to \R^{n}$ be a minimal surface and $ds^{2}$ the metric on $M$ induced from $\R^{n}$. 
According to (\ref{curvature2}), the Gauss curvature $K_{ds^{2}}$ is nonpositive. 
\begin{definition}
The {\em total curvature} of a minimal surface $M$ is defined by
\[
\tau(M)=\displaystyle \int_{M} K_{ds^{2}} dA \quad (\geq -\infty)
\]
where $dA$ is the area form of ($M$, $ds^{2}$). 
\end{definition} 
According to (\ref{metric}) and (\ref{curvature2}), the total curvature of $M$ is
\begin{equation}\label{totalcurvature}
\tau(M)=\displaystyle \int_{M} K_{ds^{2}} dA = -2\displaystyle \int_{M} \frac{|\phi \wedge {\phi}'|^{2}}{|\phi|^{4}}du\,dv\:.
\end{equation}
The relation between the Gauss map and the total curvature can be made explicit by introducing 
the Fubini-Study metric on $\Q^{n-2}(\C)$. On $\mathbb{P}^{n-1}(\C)$, 
there is a unique unitary invariant K\"{a}hler metric called the Fubini-Study metric which can be written as 
\[
{ds}^{2}=2\frac{|w\wedge dw|^{2}}{|w|^{4}}
\] 
in homogeneous coordinates $(w^{1}:\cdots:w^{n})$ of $\mathbb{P}^{n-1}(\C)$, where 
$|w\wedge dw|^{2}=\sum_{i<j}|w^{i}dw^{j}-w^{j}dw^{i}|^{2}$. The induced metric on the image under the 
Gauss map $g$ is 
\[
{ds}^{2}=2\frac{|\phi \wedge {\phi}'|^{2}}{|\phi|^{4}}|dz|^{2}\:.
\]
If we denote the image area under the Gauss map as $A(g)$, then 
\begin{equation}\label{Fubini}
A(g)=-\tau(M)\:.
\end{equation}
When the total curvature of a complete minimal surface is finite, the surface is called {\em an algebraic minimal surface}.

Now, we give the following characterization of algebraic minimal surfaces. 
\begin{theorem}[Chern-Osserman\cite{CO}, Huber\cite{H}]\label{Chern-Osserman}
An algebraic minimal surface $x\colon M\to \R^{n}$ satisfies the followings :
\begin{enumerate}
\item[(i)] $M$ is conformally equivalent to $\overline{M}\backslash\{p_{1},\ldots,p_{k}\}$ where 
$\overline{M}$ is a compact Riemann surface, and $p_{1},\ldots,p_{k}$ are finitely many points of $\overline{M}$. 
\item[(ii)] Each $\phi_{i}$ $(1\leq i\leq n)$ can be extended to $\overline{M}$ as a meromorphic $1$-form.
\end{enumerate}
\end{theorem}
\begin{proposition}[Chern-Osserman, \cite{CO}]\label{order2}
Let $x\colon M=\overline{M}\backslash\{p_{1},\ldots,p_{k}\}\to\R^{n}$ be an algebraic minimal surface. 
Then $\phi=(\phi_{1},\ldots,\phi_{n})$ has a 
pole of order $\mu_{j} \geq 2$ at each $p_{j}$.
\end{proposition}
\begin{proof}
Let $D_{j}=\{z\in\C\,|\,|z|< 1\}$ be the local coordinate centered at $p_{j}\in \overline{M}$. 
At $z=0$, the functions $f_{i}$ have at worst a pole of order $\nu_{ij}$, and so we can write 
\[
\sum|f_{i}|^{2}=\frac{c}{|z|^{2\nu_{j}}} + \mathrm{higher\;order\;terms}
\]
for some $c>0$ and $\nu_{j}=\max\{\nu_{1j},\ldots,\nu_{nj}\}\geq 1$. Assume that $\nu_{j}=1$. 
Then for some constants $(c_{1},\cdots,c_{n})\in \C^{n}$, 
\[
f_{i}=\frac{c_{i}}{z}+\mathrm{higher\;order\;terms}\quad(i=1,\ldots,n)\;.
\]
Note $\sum{f_{i}}^{2}=0$ forces
\begin{equation}\label{constant}
\displaystyle \sum_{i=1}^{n}(c_{i})^{2}=0
\end{equation}
We put
\[
\psi_{i}=f_{i}-\frac{c_{i}}{z}
\]
so that each $\psi_{i}$ is a holomorphic function near $z=0$. Note that 
\[
\Re(c_{i}\log z)=\Re\int (f_{i}-\psi_{i})dz=x_{i}-\Re\int\psi_{i}dz\;.
\]
Hence in a punctured neighborhood of $z=0$, the real part of the function $c_{i}(\log z)$ is a well-defined 
harmonic function. However imaginary part of the complex logarithm is a multi-valued function near the origin, each $c_{i}$ 
must be real. By (\ref{constant}), each $c_{i}$ must be zero, which this is a contradiction.
\end{proof}

In this thesis, we study the Gauss map of the follwing class of complete minimal surfaces 
that includes algebraic minimal surfaces. 
\begin{definition}
We call a complete minimal surface in $\R^{n}$ {\em pseudo-algebraic}, if the following conditions 
are satisfied:
\begin{enumerate}
\item[(i)] Each  $\phi_{i}$ $(i=1,\ldots,n)$ is defined on a puncutured Riemann surface 
$M=\overline{M}\backslash \{p_{1}, \ldots, p_{k}\}$, $p_{j}\in\overline{M}$, where $\overline{M}$ is a 
compact Riemann surface.
\item[(ii)] Each $\phi_{i}$ can be extended to $\overline{M}$ as a meromorphic $1$-form.
\end{enumerate}
We call $M$ {\em the basic domain} of the pseudo-algebraic minimal surface under consideration.
\end{definition}  

Since we do not assume the period condition on $M$, a pseudo-algebraic minimal surface is defined on 
some covering surface of $M$, in the worst case, on the universal covering. 
\section{The Gauss map of pseudo-algebraic minimal surfaces in $\R^{3}$}
In this section, we shall study the Gauss map of pseudo-algebraic minimal surfaces in $\R^{3}$.
\subsection{Pseudo-algebraic minimal surfaces in $\R^{3}$}
First we shall study the Gauss map of a surface in $\R^{3}$. 
In $\R^{3}$, each oriented plane $P$ is uniquely determined by the unit normal vector such that it is perpendicular 
to $P$ and the system \{X, Y, N\} is a positive orthonormal basis of $\R^{3}$ for arbitrarily chosen positively 
oriented orthonormal basis \{X, Y\} of $P$. For an oriented surface in $\R^{3}$ the tangent plane is uniquely 
determined by positively oriented unit normal vector. On the other hand, the sphere $S^{2}$ of 
all unit normal vectors in $\R^{3}$ is identified with the extended complex plane $\hat{\C}=\C\cup \{\infty\}$ 
by the stereographic projection $\varpi$.
\begin{definition}\label{CGM}
For minimal surface $M$ immersed in $\R^{3}$, the {\em Gauss map} $g:M\to \hat{\C}$ of $M$ is defined as the map 
which maps each point $p\in M$ to the point $\varpi(N_{p})\in\hat{\C}$, where $N_{p}$ is the positively oriented 
normal vector $N_{p}$ of $M$ at $p$.
\end{definition}
Next, Definition \ref{CGM} is canonically identified with the special case of that for Definition \ref{GGM}. 
Indeed, their relationship is explained as follows.
We take an arbitrary point $(w_{1}: w_{2}: w_{3})\in \Q^{1}(\C)$. Set $w_{i}=x_{i}-\sqrt{-1}y_{i}\: (1\leq i\leq 3)$ 
with real numbers $x_{i}, y_{i}$ and 
\[
W=(w^{1}, w^{2}, w^{3}),\:  X=(x^{1}, x^{2}, x^{3}),\:  Y=(y^{1}, y^{2}, y^{3})
\]
Since $(w^{1})^{2}+(w^{2})^{2}+(w^{3})^{2}=0$, they satisfy the condition (\ref{basis}). 
Multiplying $W$ by a suitable constant, we may assume that $|X|=|Y|=1$. 
Then, the unit normal vector of the plane which has a positive basis $\{X, Y\}$ is given by 
\[
N=X\times Y=\Im (w^{2}\overline{w^{3}}, w^{3}\overline{w^{1}}, w^{1}\overline{w^{2}})\;.
\]
For the case where $w^{1}\not= \sqrt{-1}w^{2}$, we assign to $W$ the point 
\[
z=\frac{w^{3}}{w^{1}-\sqrt{-1}w^{2}}
\]
and otherwise, the point $z=\infty$. This correspondence is continuous inclusively at $\infty$. 
Since $|w^{1}|^{2}+|w^{2}|^{2}+|w^{3}|^{2}=|X|^{2}+|Y|^{2}=2$, we have 
\[
\frac{1}{2}\Bigl(\frac{1}{z}-z\Bigr)=\frac{w^{1}}{w^{3}},\:  \frac{\sqrt{-1}}{2}\Bigl(\frac{1}{z}+z\Bigr)=\frac{w^{2}}{w^{3}}, \:
|w^{3}|^{2}=\frac{4|z|^{2}}{(|z|^{2}+1)^{2}}\;.
\]
These yields
\[
N=\Bigl(\frac{2\Re z}{|z|^{2}+1}, \frac{2\Im z}{|z|^{2}+1}, \frac{|z|^{2}-1}{|z|^{2}+1}\Bigr)\;.
\]
This shows that the point $S^{2}$ corresponding to $W\in\Q^{1}(\C)$ is mapped to the above point $z$ 
by the stereographic projection.

Let $x=(x^{1}, x^{2}, x^{3})\colon M\to \R^{3}$ be a surface whose the Gauss map $g$ in the sense of Definition \ref{GGM} is 
not a constant map. For a complex local coordinate $z=u+\sqrt{-1}v$, $g$ is reperesented as 
$g=(\phi_{1}: \phi_{2}: \phi_{3})=(f_{1}: f_{2}: f_{3})$, where 
\begin{equation}\label{phiphi}
\phi_{i}=f_{i}dz\; (1\leq i\leq 3)\;. 
\end{equation}
Set 
\begin{equation}\label{w-data}
hdz=\phi_{1}-\sqrt{-1}\phi_{2},\quad g=\frac{\phi_{3}}{\phi_{1}-\sqrt{-1}\phi_{2}}\;.
\end{equation}
Then the above ``$g$" is the Gauss map in the sense of Definition \ref{CGM}. 
Therefore, we see that these definitions are biholomorphically the same.
From here on we identify these two definitions. Theorem \ref{minimal2} implies :
\begin{proposition}\label{minimal23}
For a surface $M$ immersed $\R^{3}$, $M$ is a minimal surface if and only if the Gauss map is meromorphic 
on $M$.
\end{proposition}

We explain here Enneper-Weierstrass representation for minimal surfaces in $\R^{3}$. 
\begin{theorem}\label{EWformula} 
Let $x=(x^{1}, x^{2}, x^{3})\colon M\to\R^{3}$ be a non-flat minimal surface immersed in $\R^{3}$. 
Consider the holomorphic 1-forms $\phi_{1}$, $\phi_{2}$, $\phi_{3}$ and $hdz$, and the meromorphic function $g$ 
which is defined by (\ref{phiphi}) and (\ref{w-data}) respectively. Then, 
\begin{enumerate}
\item[(i)] we have
\begin{equation}\label{phis}
\phi_{1}=\frac{1}{2}(1-g^{2})hdz,\: \phi_{2}=\frac{\sqrt{-1}}{2}(1+g^{2})hdz,\: \phi_{3}=ghdz
\end{equation}
and we recover the immersion $x$ by the real Abel-Jacobi map
\begin{equation}\label{immersion}
x(z)=\Re \int^{z}_{z_{0}}(\phi_{1}, \phi_{2}, \phi_{3}) 
\end{equation}
up to translation.
\item[(ii)] the metric induced from the standard metric on $\R^{3}$ is given by 
\begin{equation}\label{metric3}
ds^{2}=\frac{|h|^{2}(1+|g|^{2})^{2}}{4}|dz|^{2}\;.
\end{equation}
\item[(iii)] the poles of $g$ of order $k$ coincides exactly with the zeros of $hdz$ of order $2k$ 
(We call it `` the regularity condition'') .
\end{enumerate}
\end{theorem}
\begin{remark}
We call the above $(hdz, g)$ {\em the Weierstrass data} of $M$.
\end{remark}

Next, we can show the following restatement of Theorem \ref{minimal3} and Theorem \ref{Chern-Osserman} 
for $n=3$. 
\begin{theorem}\label{EWformula2}
Let $M$ be an open Riemann surface, $hdz$ a non-zero holomorphic 1-form and $g$ a non-constant meromorphic function 
$M$. Assume that the poles of $g$ of order $k$ coincides exactly with the zeros of $hdz$ of order $2k$ and 
that the holomorphic 1-forms $\phi_{1}$, $\phi_{2}$, $\phi_{3}$ defined by (\ref{phis}) have no real periods. 
Then, for the functions $x^{1}$, $x^{2}$, $x^{3}$ defined by (\ref{immersion}), the surface
\[
x=(x^{1}, x^{2}, x^{3})\colon M\to\R^{3}
\] 
is a minimal surface immersed in $\R^{3}$ whose Gauss map is the map $g$ and whose induced metric is given (\ref{metric3}).
\end{theorem}

\begin{theorem}\label{Huber-Osserman}
An algebraic minimal surface
$x:M\to \R^3$ satisfies the following : 
\begin{enumerate}
\item[(i)] $M$ is conformally equivalent to 
$\overline{M}\setminus\{p_1,\dots,p_k\}$ where $\overline{M}$ 
is a compact Riemann surface, and $p_1,\dots,p_k$ are finitely many points of $\overline{M}$.
\item[(ii)]  The Weierstrass data $(hdz,g)$ is extended 
meromorphically to $\overline{M}$.
\end{enumerate}
\end{theorem}

Here, we give some examples. We denote the number of exceptional values of $g$ by $D_{g}$. 
\begin{example}[Enneper surface]
On $M=\C$, we consider the Weierstrass data 
\[
(hdz, g)=(dz, z)\;
\]
The regularity condtion satisfied and the period condition is satisfied since $M$ is simply connected. 
The resulting minimal surface $x\colon \C\to\R^{3}$ 
is called Enneper surface. It is an algebraic minimal surface with total curvature $-4\pi$ and $D_{g}=1$.
\end{example}
\begin{example}[Catenoid]
On $M=\C\backslash \{0\}$, we consider the Weierstrass data 
\[
(hdz, g)=\Bigl(\frac{dz}{z^{2}}, z\Bigr)\;.
\]
The data satisfies the regularity condtion and the period condition.
The resulting minimal surface $x\colon \C\backslash \{0\}\to\R^{3}$ 
is called catenoid. It is an algebraic minimal surface with total curvature $-4\pi$ and $D_{g}=2$.
\end{example}

Other than Catenoid, there are many examples of algebraic minimal 
surfaces with $D_g=2$,  which include those of hyperbolic type. 

\begin{theorem}[Miyaoka-Sato \cite{MS}]
There exist algebraic minimal surfaces with $D_g=2$, for
\begin{enumerate}
\item[\rm{(i)}] $G=0, k\geq 2$
\item[\rm{(ii)}]  $G=1, k\geq 3$
\item[\rm{(iii)}]  $G\geq 2, k\geq 4$
\end{enumerate}
\end{theorem}

When $G=0$ and $k=2$, all such minimal surfaces are classified.
Examples for $G=0$ and $k=3$  
given below [\cite{MS}, Proposition 3.1] are important for later argument:
 let $M=\textbf{P}^1\setminus\{\pm i,\infty\}$,
and define a Weierstrass data by
\begin{equation}\label{Miyaoka-Sato}
\left\{
\begin{array}{ll}
g(z)=\sigma\dfrac{z^2+1+a(t-1)}{z^2+t}\\
hdz=\dfrac{(z^2+t)^2}{(z^2+1)^2}dz,\quad (a-1)(t-1)\ne 0\\
\sigma^2=\dfrac{t+3}{a\{(t-1)a+4\}}\,.
\end{array}
\right.
\end{equation}
For any $a,t$ satisfying $\sigma^2<0$, we obtain an algebraic
minimal surface whose Gauss map omits two values 
$\sigma,\sigma a$.

Applying the covering method to this surface 
(see Remark 3.21), we obtain examples 
of (ii) and (iii).
However as these examples have all the same image in $\R^3$, Miyaoka and Sato further 
 constructed mutually non-congruent examples for $G=1$ and $k=4$, by generalizing Costa's surface 
[\cite{MS},Theorem 3].  For details see Remark 3.22.

\begin{example}[Helicoid]
On $M=\C$, we consider the Weierstrass data 
\[
(hdz, g)=(e^{-z}dz, \sqrt{-1}e^z) \;.
\]
The regularity condition satisfied and the period condition is vacuously satisfied. 
The resulting minimal surface $x\colon \C\to\R^{3}$ 
is helicoid. It is a minimal surface with infinite total curvature and $D_{g}=2$.
\end{example}
\begin{example}[Jorge-Meeks surface\cite{JM}]
Let $\Sigma_{r}=\{z\in\C |z^{r}=1\}$. We take $M=\hat{\C}\backslash \Sigma_{r}$ and consider 
\[
(hdz, g)=\Bigl(\frac{dz}{(z^{r}-1)^{2}}, z^{r-1}\Bigr) \;.
\]
Then we can show that $(hdz, g)$ is a Weierstrass data defining an algebraic minimal surface with total curvature 
$-4(r-1)\pi$ and when $r\geq 3$, $D_{g}=0$. 
\end{example}

\begin{example}[Costa surface]\label{costa}
Let $\overline{M}$ be the square torus on which the Weierstrass $\wp$ functions satisfies 
$(\wp')^{2} = 4\wp({\wp}^{2}-a^{2})$. Let $M$ be given by removing $3$ points satisfying $\wp=0, \pm a$ from 
$\overline{M}$. On $M$, we consider the Weierstrass data
\[
(hdz, g)=\Bigl(\wp(z)dz, \frac{\wp(z)}{\wp'(z)}\Bigr)\;,
\] 
where $A=2\sqrt{2\pi}\wp(1/2)$. Then we can show that the minimal surface defined by this Weierstrass data is 
an algebraic minimal surface with total curvature $-12\pi$ and $D_{g}=1$.
\end{example}

In general, for a given meromorphic function $g$ on $M$, it is not so hard to find a holomorphic 1-form $hdz$ 
satisfying the regularity condition. However, the period condition always causes trouble. When the period condtion is not 
satisfied, we anyway obtain a minimal surface on the universal covering surface of $M$. 

Here we notice that the triple of holomorphic $1$-forms 
\[
e^{i\theta}(\phi_1,\phi_2,\phi_3),\quad \theta\in \R
\]
also satisfies the regularity condition.
The corresponding  Weierstrass data is given by
\begin{equation}\label{associate}
\left\{
\begin{array}{ll}
g^\theta(z)=g(z)\\
h^\theta dz=e^{i\theta}hdz\,.
\end{array}
\right.
\end{equation}
As the period condition is scarcely satisfied by these data,
we get an  $S^1$ parameter family of minimal surfaces defined 
on the universal covering surface by (\ref{immersion}), which is called 
the associated family.
Note that all surfaces in this family have the same 
Gauss map.

As easily seen from (\ref{curvature2}) and the assertion (ii) of Theorem \ref{EWformula}, the Gauss curvature 
of $M$ is given by 
\begin{equation}\label{gausscurv}
K_{ds^{2}}(p)=-\frac{4|g'|^{2}}{|h|^{2}(1+|g|^{2})^{4}}
\end{equation}
and the total curvature by 
\begin{equation}\label{totalcurv3}
\tau(M)=-\int_{M}\Bigl(\frac{2|g'|}{1+|g|^{2}}\Bigr)^{2} du\wedge dv =-4\pi d,\; d\in \N\cup \{\infty\}\;.
\end{equation}

For a complete minimal surface in $\R^{3}$, the definition of ``pseudo-algebraic'' is as follows.
\begin{definition}
We call a complete minimal surface in $\R^3$ {\em pseudo-algebraic}, 
if the following conditions are satisfied:
\begin{enumerate}
\item[(i)] The Weierstrass
data $(hdz,g)$ is defined on a Riemann surface $M=\overline{M}\setminus\{p_1,\dots,p_k\}$, $p_j\in \overline{M}$, 
where $\overline{M}$ is a compact Riemann surface.
\item[(ii)] $(hdz,g)$ can be extended meromorphically to $\overline{M}$.
\end{enumerate}
We call $M$ {\em the basic domain} of the pseudo-algebraic minimal surface under consideration. 
\end{definition}

\begin{remark}
Gackst\"atter called such surfaces {\em abelian minimal surfaces} \cite{G}.
\end{remark}

Algebraic minimal surfaces and their associated surfaces are
certainly pseudo-algebraic. Another important example is Voss' surface.
The Weierstrass data of this surface is defined on $M=\C\setminus\{a_1,a_2,a_3\}$ for
distinct $a_1,a_2,a_3\in \C$, by
\begin{equation}\label{Voss}
\left\{
\begin{array}{ll}
g(z)=z\\
hdz=\dfrac{dz}{\Pi_j(z-a_j)}\,.
\end{array}
\right.
\end{equation}
As this data does not satisfy the period condition,
we get a minimal surface $x:\D\to \R^3$ on the universal
covering disk of $M$.
In particular, it has infinite total curvature.
We can see that the surface is complete
 and the Gauss map omits four values $a_1,a_2,a_3,\infty$.
Starting from $M=\C\setminus\{a_1,a_2\}$, we get similarly 
a complete minimal surface $x:\D\to \R^3$, 
of which Gauss map omits three values $a_1,a_2,\infty$.
The completeness restricts the number of points $a_j$'s 
to be less than four. 
In  both cases, all elements in the associated 
family have infinite total curvature.
\begin{remark} 
There exists a complete minimal surface which is ``not" pseudo-algebraic with $D_{g}=4$. 
For details see \cite{L}.
\end{remark}

\subsection{Ramification estimate and unictiy theorem}
\begin{definition}
We call $b\in \mathbb{P}^1(\C)$ a totally ramified value of $g$ when at any 
inverse image of $b$, $g$ branches. 
We regard exceptional values also as totally ramified values.
Let $\{a_1,\dots,a_{r_o},b_1,\dots,b_{l_0}\}\subset \mathbb{P}^1(\C)$ be the set of 
 totally ramified values of $g$, where $a_j$'s are exceptional 
values. For each $a_j$, put $\nu_j=\infty$, and for each $b_j$, 
define $\nu_j$ to be the minimum of the multiplicity of $g$ at 
points $g^{-1}(b_j)$. Then  we have $\nu_j\geq 2$.
We call 
\[
\nu_{g}=\sum_{a_j,b_j}(1-\dfrac1{\nu_j})=r_0+\sum_{j=1}^{l_0}(1-\dfrac1{\nu_j})
\]
{\em the totally ramified value number of $g$}.
\end{definition}

A natural meaning of this number is explained in the framework of 
we need the second main theorem in the Nevanlinna theory.
We refer to \cite{Ko} for this theory. 
Note that though $\nu_{g}$ is a rational number. Fujimoto proved the following.
\begin{theorem}[Fujimoto \cite{F4}] \label{ino} 
Let $x\colon M\to\R^{3}$ be a non-flat complete minimal surface, $g$ be its Gauss map. 
Then we have
\[
D_{g}\leq \nu_{g}\leq 4\;.
\]
\end{theorem}
On the other hand, Osserman proved the following theorem.
\begin{theorem}[Osserman \cite{O1}]
Let $x\colon M\to\R^{3}$ be a non-flat algebraic minimal surface, $g$ be its Gauss map. 
Then we have
\[
D_{g}\leq 3
\]
\end{theorem}
However, there is no known example with $D_{g}=3$. Since there are many example with $D_{g}=2$, 
many people believe ``2" is the best possible upper bound of $D_{g}$. Moreover, as in the case of Fujimoto's 
theorem (Theorem \ref{ino}), it has been implicitly believed that the same is true for $\nu_{g}$. 
However, we discovered that this is false by the following result. 
\begin{theorem}[Kawakami \cite{Ka}]
The Gauss map of the algebraic minimal surfaces 
given in $(\ref{Miyaoka-Sato})$ 
has totally ramified value number 2.5. 
\end{theorem}

In fact, it has two 
exceptional values, and another totally ramified value at $z=0$ where 
$g'(z)=0$.

Now, we give the ramification estimates of the Gauss map of pseudo-algebraic minimal surfaces in $\R^{3}$.
\begin{theorem}[Kawakami, Kobayashi and Miyaoka \cite{KKM}] \label{Thm3.1} 
Consider a pseudo-algebraic minimal surface with the basic domain $M=\overline{M}\setminus\{p_1,\dots,p_k\}$.
Let $G$ be 
the genus of $\overline{M}$, and let $d$ be the degree of $g$
considered as a map on $\overline{M}$.
Then we have
\begin{equation} \label{Exp1}
D_g\leq 2+\dfrac2{R}\,,\quad
R=\dfrac{d}{G-1+k/2}\geq 1\,.
\end{equation}
More precisely, if the number of (not necessarily totally) 
ramified values other than the exceptional values of $g$ is $l$, 
we have 
\begin{align}\label{Exp2}
D_g\leq 2+\dfrac2{R}-\dfrac{l}{d}\,.
\end{align}
On the other hand, the totally ramified value number of $g$  
satisfies 
\begin{equation}\label{TRVN1}
 \nu_{g}\leq 2+\dfrac2{R}\,. 
\end{equation}
In particular, we have
\begin{equation}\label{Total1}
D_g\leq \nu_{g}\leq 4\,,
\end{equation}
and for algebraic minimal surfaces, the second inequality is a 
strict inequality. 
(\ref{Exp2}) and (\ref{TRVN1}) are best possible in both 
algebraic and non-algebraic cases. 
\end{theorem}
\begin{proof}
The proof of Theorem \ref{Thm3.1} is given by a refinement of the proof 
of Osserman's theorem in \cite{O1}. In order to simplify the argument, 
we may assume without loss of generality 
that $g$ is neither zero nor pole at  $p_j$, and moreover, 
zeros and poles of $g$ are simple.
By completeness, $hdz$ has poles of order $\mu_j\ge 1$ at $p_j$.
By Proposition \ref{order2}, the period condition implies $\mu_j\geq 2$, 
however here we do not assume this.
Let $\alpha_s$ be (simple) zeros of $g$, $\beta_t$ (simple) poles of $g$.
The following table shows the relation between zeros and poles 
of $g$, $hdz$ and $ghdz$. 
The upper index means the order. 
\begin{center}
\begin{tabular}{|c|c|c|c|}\hline
$z$ & $\alpha_s$  & $\beta_t$ & $p_j$  \\\hline
$g$& $0^1$ & $\infty^1$ &   \\\hline
$hdz$ &   & $0^2$ & $\infty^{\mu_j}$  \\\hline
$ghdz$ & $0^1$ & $0^1$ & $\infty^{\mu_j}$ \\\hline
\end{tabular}
\end{center}
Applying the Riemann-Roch formula to the meromorphic
differential $hdz$ or $ghdz$ on $\overline{M}$, we obtain
\[
2d-\sum_{j=1}^k\mu_j=2G-2\,.
\]
Note that this equality depends on the above setting of
zeros and poles of $g$, though $d$ is an invariant.
Thus we get
\begin{equation}\label{Riemann-Roch}
d=G-1+\dfrac12\sum_{j=1}^k\mu_j\ge G-1+\dfrac{k}2\,,
\end{equation}
and
\begin{equation}\label{ration}
R\ge 1\,.
\end{equation}
When $M$ is an algebraic minimal surface or its associated surface,
we have $\mu_j\ge 2$ and so $R>1$.

Now, we prove (\ref{Exp2}) (and (\ref{Exp1})).
Assume $g$ omits $r_0=D_g$ values, and let $n_0$ be the sum of the 
branching orders of $g$ at these exceptional values.
Moreover, let $n_b$ be the sum of branching orders at the inverse
images of  
non-exceptional (not necessarily totally) ramified values $b_1,\dots,b_l$ of $g$.  
We see
\begin{equation}\label{exbr}
k\ge dr_0-n_0\,,\quad 
n_b\ge l\,.
\end{equation}
Let $n_{g}$ be the total branching order of $g$. 
Then applying Riemann-Hurwitz's theorem to the 
meromorphic function $g$ on $\overline{M}$,
we obtain
\begin{equation}\label{Riemann-Hurwitz}
n_{g}=2(d+G-1)=n_0+n_b\ge dr_0-k+l\,.
\end{equation}
If we denote
\[
\nu_i=\text{min}_{g^{-1}(b_i)}\{\text{multiplicity of }g(z)=b_i\}\,,
\]
we have $1\le \nu_i\le d$. 
Now the number of exceptional values satisfies
\begin{equation}\label{Exp3}
D_g= r_0\leq \dfrac{n_{g}+k-l}{d}=2+\dfrac2{R}-\dfrac{l}{d}
\end{equation}
where we have used (\ref{Riemann-Hurwitz}), hence (\ref{ration}) implies
\[
D_g\leq 2+\dfrac2{R}\leq 4\,.
\]
In particular for algebraic minimal surfaces and its associated 
family, we have $R>1$ so that
\[
D_g\le 3\,,
\]
which is nothing but Osserman's theorem.

Next, we show (\ref{TRVN1}).
Let $b_1,\dots,b_{l_0}$ be the {\em totally} ramified values which are not exceptional values. Let $n_r$ be the sum of branching orders 
at $b_1,\dots,b_{l_0}$. For each $b_i$, the number of points
in the inverse image $g^{-1}(b_i)$ 
is less than or equal to $d/\nu_i$, since $\nu_i$ is the 
minimum of the multiplicity at all $g^{-1}(b_i)$.
Thus we obtain
\begin{equation}\label{trvn}
dl_0-n_r\leq \sum_{i=1}^{l_0}\dfrac{d}{\nu_i}\,.
\end{equation}
This implies
\[
l_0-\sum_{i=1}^{l_0}\dfrac{1}{\nu_i}\leq \dfrac{n_r}{d}\,,
\]
hence using the first inequality in (\ref{exbr}) and $n_r\leq n_b$, we get
\[
\nu_{g}=r_0+\sum_{i=1}^{l_0}(1-\dfrac{1}{\nu_i})\leq 
\dfrac{k+n_0}{d}+\dfrac{n_r}{d}\leq \dfrac{n_{g}+k}{d}=2+\dfrac2{R}\,.
\]
\end{proof}

The sharpness of (\ref{Exp2}) and (\ref{TRVN1}) follows from:
\begin{enumerate}
\item[(1)] When $d=2$ we have
\[
D_g\le 2+\frac2{R}-\frac{l}{2}, \quad \nu_{g}\leq 2+\dfrac2{R}\,.
\]
The surface given by (\ref{Miyaoka-Sato}) attains both equalities, since $R=4$, $l=1$ and
$D_g=2$, $\nu_{g}=2.5$. Thus (\ref{Exp2}) and (\ref{TRVN1}) are sharp.
\item[(2)] 
Voss' surface satisfies $d=1$ and $G=0$.
Thus when $k=3$, we get $R=2$, $l=0$ hence 
$D_g=3=2+2/2$.
When $k=4$, we have $R=1$, $l=0$ and 
$D_g=4=2+2/1$. These show that (\ref{Exp2}) and (\ref{TRVN1}) are
sharp in non-algebraic pseudo-algebraic case, too.
\end{enumerate}

\begin{remark} \label{coveringmethod}
There exists a way of construction of algebraic minimal surfaces 
by a covering method of Klotz-Sario \cite{BC}.
Indeed, if $x:M\to \R^3$ is an algebraic minimal surface,
and if $\pi:\hat M\to M$ is a non-branched covering surface 
of $M=\overline{M}\setminus\{p_1,\dots,p_k\}$,
then we obtain a new algebraic minimal surface by 
$\hat x=x\circ \pi:\hat M\to \R^3$. 
This surface has the same image as the original one, but
the domain $\hat M$ has different topological type.
Nevertheless, we can see that {\em the ratio $R$ is invariant}
under this construction, via a little algebraic argument.
Certainly, $D_g$ and $\nu_{g}$ are also invariant under covering
construction.
\end{remark}

\begin{remark}
The inequality (\ref{Exp1}) is also best possible 
for algebraic minimal surfaces in the following sence. 
In [\cite{MS}, Theorem 3], Miyaoka and Sato constructed two infinite 
series of mutually distinct algebraic minimal surfaces 
of the fixed topological type $G=1$ and $k=4$, 
whose Gauss map omits 2 values.  
These surfaces are given as follows. Let $\overline{M}$ be the square torus  
on which the Weierstrass $\wp$ 
function satisfies $(\wp')^2=4\wp({\wp}^2-a^2)$. 
Let $M$ be given by removing 
4 points satisfying $\wp=0,\pm a,\infty$ from $\overline{M}$.
Define the Weierstrass data by \\
(Case 1) $g=\dfrac\sigma{{\wp}^j{\wp}'},\quad 
hdz=\dfrac{{\wp}d{\wp}}{{\wp}'},\quad j=1,2,3\dots,$\\
(Case 2) $g=\dfrac\sigma{{\wp}^j{\wp}'},\quad 
hdz=\dfrac{{\wp}^{j+1}d{\wp}}{{\wp}'},\quad j=2,4,6\dots,$ \\
Then choosing a suitable $\sigma$, we obtain algebraic minimal 
surfaces with $g$ omitting 2 values 0 and $\infty$.
Since the degree of $g$ is $d=2j+3$ in both cases and  
$R=d/2=(2j+3)/2$,  
$2+2/R$ tends to $2$ $(=D_g)$ as close as we like. 
(Costa's surface (Example \ref{costa}) is given by $j=0$, in which case
$(G,k,d)=(1,3,3)$, and $g$ omits just one value 0.)
\end{remark}

\begin{remark} 
The inequality (\ref{Exp2})  gives
more informations than (\ref{Exp1}). 
In particular, (\ref{Exp2}) implies that the more branch points 
$g$ has in $M$,  the less is the number of exceptional 
values. 
\end{remark}

\begin{remark}
When we prove Throrem \ref{Thm3.1} for algebraic minimal surfaces, 
we use ``local'' period conditions as $\mu_{j}\geq 2$. 
However we do not use ``global'' period conditions i.e., the element of $H_{1}(\overline{M},\Z)$. 
Thus we do not understand how this affect the estimete of $D_{g}$ or $\nu_{g}$.
This is our future problem.
\end{remark}
The geometrical meaning of the ration ``$R$'' is given in the next subsection. 
Theorem \ref{Thm3.1} implies the following known facts:
\begin{corollary} [cf.\:Osserman \cite{O1}, Fang \cite{Fa}, Gackst\"atter \cite{G}] \label{threepoint}
For algebraic minimal surfaces, we have:
\begin{enumerate}
\item[(i)] When $G=0$, $D_g\le 2$ holds.
\item[(ii)] When $G=1$ and $M$ has a non-embedded end, $D_g\le 2$ holds.
If $G=1$ and $D_g=3$ occur, $d=k$ follows and $g$ does not branch 
in $M$, so is a non-branched covering of $\mathbb{P}^1(\C)\setminus \{3\text{ points}\}$.
\end{enumerate}
\end{corollary}
\begin{proof}
 It is easy to see that $r_0=3$ implies $R\leq 2$, hence 
\[
G-1+\dfrac12\sum_{j=1}^k\mu_j\leq 2(G-1)+k\; .
\]
As we have $\mu_j\geq 2$ in the algebraic case, it follows
\begin{equation}\label{O-ineq.}
k\leq \dfrac12\sum_{j=1}^k\mu_j\leq G-1+k\,.
\end{equation}
Thus  we obtain (i). 
When $G=1$, we get $\mu_j=2$ for all $j$, which means that 
all the ends are embedded (\cite{JM}), and $R=2$. But since
$R=\dfrac{d}{k/2}$, we obtain $d=k$. Finally from 
(\ref{Exp2}), we get $l=0$, which means that $g$ does not
branch in $M$.
\end{proof}
\begin{remark}
Fang [\cite{Fa}, Theorem 3.1] shows that algebraic 
minimal surfaces with $d\leq 4$ satisfy $D_g\leq 2$  (see \cite{WX}
for $d\leq 3$).
\end{remark}

Here, We give two applications of Theorem \ref{Thm3.1}.
First one is a unicity theorem for the Gauss map of pseudo-algebraic minimal surfaces.
\begin{theorem}[Kawakami, Kobayashi and Miyaoka \cite{KKM}]
Consider two pseudo-algebraic minimal surfaces $M_1, M_2$ with the 
same basic domain $M=\overline{M}\setminus\{p_1,\dots,p_k\}$.
Let $G$ be the genus of $\overline{M}$, and let 
$g_1, g_2$ be the Gauss map of $M_1$ and $M_2$ respectively.
Assume that $g_1$ and $g_2$ have the same degree $d$ 
when considered as a map on $\overline{M}$, but assume $g_1\ne g_2$ 
as a map $M\to \mathbb{P}^1(\C)$. Let $c_1,\dots, c_q\in \mathbb{P}^1(\C)$ be
distinct points such that $g_1^{-1}(c_j)\cap M=g_2^{-1}(c_j)\cap M$ for 
$1\leq j\leq q$. Then 
\begin{equation}\label{unicity1}
q\leq 4+\frac2{R},\quad R=\dfrac{d}{G-1+k/2}
\end{equation}
follows. In particular, $q\leq 6$, and for algebraic 
minimal surfaces we have $q\leq 5$.
\end{theorem}
\begin{proof}
Put
\[
\delta_j=\sharp(g_1^{-1}(c_j)\cap M)=\sharp (g_2^{-1}(c_j)\cap M)\,,
\]
where $\sharp$ denotes the number of points.
Then we have
\begin{equation}\label{branchcount}
qd\leq k+\sum_{j=1}^q\delta_j+n_{g}\,,
\end{equation}
using the same notation as in proof of Theorem \ref{Thm3.1}.
Consider a meromorphic function $\varphi=\dfrac1{g_1-g_2}$ on $M$.
Then at each point of $g_1^{-1}(c_j)\cap M$, $\varphi$ has 
a pole, while the total number of the poles of $\varphi$ 
is at most $2d$, hence we get
\begin{equation}\label{function}
\sum_{j=1}^q\delta_j\leq 2d\,.
\end{equation}
Then from (\ref{branchcount}) and (\ref{function}), we obtain
\[
qd\leq k+2d+n_{g}\,,
\]
and
\[
q\leq \dfrac{2d+n_{g}+k}{d}=4+\dfrac2{R}
\]
follows immediately.
\end{proof}
\begin{remark}
Fujimoto \cite{F3} gives an example of two 
pseudo-algebraic minimal surfaces with $q=6$, of which Gauss maps
do not coincide. For algebraic case, whether $q=5$ is best possible or
not is an interesting open problem.
\end{remark}

Next, for later use, we mention Gackst\"atter's result \cite{G} :
\begin{proposition}[Gackst\"atter \cite{G}] 
If the Gauss map of an algebraic minimal surface with $G=1$ 
omits 3 values $a_1,a_2,a_3\in \mathbb{P}^1(\C)$, 
then all branch points of $g$ are located at the end points, and
$g$ is  a non-branched covering map 
 of $\mathbb{P}^1(\C)\setminus\{a_1,a_2,a_3\}$. 
\end{proposition}

This follows immediately from Corollary \ref{threepoint} (ii). 
Thus the Gauss map descends to 
$\mathbb{P}^1(\C)\setminus \{3\text{ points}\}$, but the minimal surface 
is not obtained from a covering of a minimal surface 
defined on $\mathbb{P}^1(\C)\setminus \{3\text{ points}\}$,
otherwise, by (ii) of Corollary \ref{threepoint}, $D_g\leq 2$. 
This implies that $hdz$ can not descends to 
$\mathbb{P}^1(\C)\setminus \{3\text{ points}\}$.

The following is obvious:
\begin{proposition}[Kawakami, Kobayashi and Miyaoka \cite{KKM}]
If the Gauss map $g$ of a pseudo-algebraic minimal surface 
omits $r$ values $a_1,\dots,a_r\in \mathbb{P}^1(\C)$ for $r=3,4$, and has 
no branch points in the basic domain $M$, then $g$ is 
a non-branched covering of $\mathbb{P}^1(\C)\setminus \{a_1,\dots,a_r\}$.
\end{proposition}

Since $r\geq 3$, the universal covering 
surface of $M$ and of $\mathbb{P}^1(\C)\setminus\{a_1,\dots,a_r\}$ are disks,
which we denote by $\D$ and $\Omega$, respectively.
When $g$ has no branch points in $M$, 
 the lifted map
$g:\D\to \Omega$ is a non-branched holomorphic map, 
i.e., a hyperbolic isometry. 
Since the degree of $g$ restricted to $\overline{M}$ is $d$, 
the fundamental domain of $M$ is given by $\cup_{j=1}^d U_i\subset \D$,  where each $U_i$ is diffeomorphic to 
$\mathbb{P}^1(\C)\setminus \{a_1,\dots,a_r\}$. 
\begin{example}
Voss' surfaces are examples for $d=1$.
\end{example}
\subsection{Some results on Nevanlinna theory of the Gauss map }
In this subsection, We state some links to the Nevanlinna theory. 
 
We consider the case where the universal covering 
surface of $M$ is a unit disk $\D$.
In order to adjust to the Nevanlinna theory, we use the 
 hyperbolic metric $\omega_h$ with curvature 
$-4\pi$ on $\D$, and 
the Fubini-Study metric $\omega_{FS}$ with curvature $4\pi$ on $\mathbb{P}^1(\C)$ 
(hence $\mathbb{P}^1(\C)$ has area 1).
Then by Gauss-Bonnet's theorem for a complete punctured Riemann
surfaces with hyperbolic metric, we have
\begin{equation}\label{GaussBonnet}
2\pi\chi(M)=\int_M K_h \omega_h=-4\pi \int_M \omega_h
=- 4\pi A_{hyp}(M)\,,
\end{equation}
where $A_{hyp}(M)$ is 
the hyperbolic area of $M$, hence for the fundamental domain $F$ of 
$M$, we get
\begin{equation}\label{hyperbolic}
A_{hyp}(F)=G-1+\dfrac{k}{2}\,.
\end{equation}

\begin{remark}
The Gauss-Bonnet theorem (\ref{GaussBonnet}) for a punctured Riemann surface $(M,\omega_h)$ is often used without proof, so here we give a brief proof. 
Let $D_{\varepsilon_j}$
be the disk with radius $\varepsilon_j$ around $p_j$, $j=1,2,\dots, k$.
We denote $M_{\varepsilon}=\overline{M}\setminus \cup_jD_{\varepsilon_j}$,
and by $\varepsilon\to 0$, we mean all $\varepsilon_j\to 0$.
Consider any metric $\sigma$ on $\overline{M}$ which is flat in all $D_{\varepsilon_j}$. Denoting locally (as K\"ahler forms) 
$\sigma=\dfrac{i}2\tilde \sigma dz\wedge d\bar z$  
and $\omega_h=\dfrac{i}2\tilde \omega_hdz\wedge d\bar z$, 
we have by Stokes' theorem 
\[
-\sum_j\int_{\partial D_{\varepsilon_j}}d^c\log (\sigma/\omega_h)
=\int_{M_{\varepsilon}}dd^c\log (\sigma/\omega_h)
=\int_{M_{\varepsilon}}dd^c\log \tilde\sigma
-\int_{M_{\varepsilon}}dd^c\log \tilde\omega_h\,,
\]
where $d=\partial+\bar\partial$, $d^c=(\partial-\bar\partial)/(4\pi i)$, 
(here $\partial$ is the half of Osserman's one). 
Because $dd^c\log\tilde\omega=-\dfrac{K_\omega}{2\pi}dA_{\omega}$ 
holds where $K_\omega$ and $dA_\omega$ are the Gauss curvature 
and the area form of $\omega$, respectively, 
taking the limit $\varepsilon\to 0$ and applying 
the Gauss-Bonnet's theorem to $(\overline{M},\sigma)$, we obtain  
\[
\lim_{\varepsilon\to 0}\sum_j
\int_{M_{\varepsilon}}(dd^c\log \tilde\sigma-dd^c\log \tilde\omega_h)
=-\chi(\overline{M})-2A_{hyp}(M)\,.
\]
Next, take a local coordinate on each $D_{\varepsilon_j}$ 
so that $z=0$ corresponds to $p_j$. Then we can express 
$\sigma=\dfrac{i}2dz\wedge d\bar z$ 
and
 $\omega_h=\dfrac{i}{2\pi}\dfrac{dz\wedge d\bar z}{|z|^2(\log|z|^{-2})^2}$ on $D_{\varepsilon_j}$.
Noting that $d^c=\dfrac1{4\pi}\bigl(-\dfrac1{r}\dfrac{\partial}{\partial\theta}dr+r\dfrac{\partial}{\partial r}d\theta\bigr)$, we obtain 
\[
\lim_{\varepsilon\to 0}\sum_j\int_{\partial D_{\varepsilon_j}} d^c\log (\sigma/\omega_h)=k\,,
\]
which implies (\ref{GaussBonnet}) and (\ref{hyperbolic}).
\qed
\end{remark}
Next, let $d$ be the degree of $g$, 
then the area $A_{FS}(F)$ of $F$
with respect to the induced metric $g^*\omega_{FS}$ is  $d$.
Thus we obtain
\begin{equation}
A_{FS}(F)=\dfrac{d}{G-1+k/2}A_{hyp}(F)=RA_{hyp}(F)\,.
\end{equation}
We now know the meaning of the ratio $R$;  the ratio of 
the area of the fundamental domain with respect to the 
induced Fubini-Study metric
to the one with respect to the hyperbolic metric on $\D$.
\begin{remark}
Even when the conformal type of $M$ is not  
hyperbolic, the ratio $R$ is meaningful in 
Theorem \ref{Thm3.1}.
\end{remark}

Now, we recall Shimizu-Ahlfors' theorem on the characteristic
function $T_g(r)$ of $g$, which states
\[
T_g(r)=\int_0^r \dfrac{dt}{t}\int_{\C(t)}g^*\omega_{FS}\,.
\]
Here $\C(t)$ is the subdisk of $\D$ 
with radius $0<t<1$. 
In order to develop the Nevanlinna theory on meromorphic functions 
on the unit disk, we need the growth order of  $T_g(r)$ 
compared with
\[
\int_0^r\dfrac{dt}{t}\int_{\C(t)}\omega_{h}\approx
\dfrac12\log\dfrac1{1-r}\,,
\]
where $r$ is sufficiently close to 1
(strictly, the left hand side is
$\dfrac12\log\dfrac1{1-r^2}$). 
We always use this approximation formula
in the following discussion, because in the Nevanlinna theory,
a bounded quantity is ignored.

If we replace the Fubini-Study 
metric by a singular metric $\Psi$ on $\mathbb{P}^1(\C)$ with area 1, 
we have
\begin{equation}
T_g(r)\geq \displaystyle\int_0^r\dfrac{dt}{t}\int_{\C(t)}g^*\Psi\,.
\end{equation}
This is shown rather easily by using Crofton's formula 
in the integral geometry \cite{Ko}.
When the image $g(M)$ is $\mathbb{P}^1(\C)\setminus \{r\text{ points}\}$,
where $r=3$ or $4$, the singular metric $\Psi$ on $\mathbb{P}^1(\C)$ 
induced by the hyperbolic metric on $\Omega$ normalized so that the 
area of $g(M)$ (counted without multiplicity) is $1$ fits the case.
Using this metric, we give a few computable examples. 
\begin{proposition}[Kawakami, Kobayashi and Miyaoka \cite{KKM}]
Consider a pseudo-algebraic minimal surface with 
the basic domain $M=\overline{M}\setminus\{p_1,\dots,p_k\}$, 
and assume that $g$ branches only at $p_j$'s.   
\begin{enumerate}\label{3point}
\item[\rm{(i)}] If $D_g=3$, we have 
\begin{equation}
T_g(r)\geq \log\dfrac1{1-r}\,.
\end{equation}
This is satisfied by Voss' surface with $k=3$, 
and an algebraic minimal surface with $G=1$ and $D_g=3$, if any.
\item[\rm{(ii)}] \label{4point}If $D_g=4$, we have 
\begin{equation}\label{character}
T_g(r)\geq \dfrac12\log\dfrac1{1-r}\,.
\end{equation}
This is satisfied by Voss' surface with $k=4$. 
\end{enumerate}
\end{proposition}
\begin{proof}
Let $\D$ be the universal covering disk of $M$, and 
$\Omega$ that of $\mathbb{P}^1(\C)\setminus\{a_1,\dots,a_{r_0}\}$, 
where $a_1,\dots,a_{r_0}$ are the exceptional values of $g$.
Let $\omega_{h}$ and $\omega_{\Omega}$ be the 
hyperbolic metric with curvature $-4\pi$.
Denote by $g:\D\to \Omega$ the lifted 
Gauss map. Since this is not branched, $g$
is a hyperbolic isometry.
To obtain the characteristic function $T_g(r)$, 
normalize the metric $\omega_{\Omega}$ so that 
the fundamental domain of $\mathbb{P}^1(\C)\setminus\{a_1,\dots,a_{r_0}\}$
has area $1$.
When $D_g=r_0=3$, this area with respect to $\omega_{\Omega}$ 
is $G-1+3/2=1/2$ by (\ref{hyperbolic}), 
thus we use the metric $2\omega_{\Omega}$ in (\ref{character}), 
and we get 
\[
\begin{array}{ll}
T_g(r)\geq\displaystyle\int_0^r\dfrac{dt}{t}\int_{\C(t)}2g^*\omega_{\Omega}
=2\displaystyle\int_0^r\dfrac{dt}{t}\int_{\C(t)}\omega_h\\
=\log\dfrac1{1-r}\,.
\end{array}
\]
The last assertion in (i) follows from Proposision 5.3.
When $D_g=r_0=4$, we need no change of the metric, and get (ii).
\end{proof}

\section{The Gauss map of pseudo-algebraic minimal surfaces in $\R^{4}$}
In this section, we shall study the Gauss map of pseudo-algebraic minimal surfaces in $\R^{4}$. 
\subsection{Pseudo-algebraic minimal surfaces in $\R^{4}$}
First we shall study the Gauss map of a surface in $\R^{4}$. The Gauss map $g$ of a surface 
in $\R^{4}$ is a holomorphic map into $\Q^{2}(\C)$. We shall inquire into the structure of $\Q^{2}(\C)$. We define 
the map $\psi_{1}\colon \Q^{2}(\C)\to \mathbb{P}^{1}(\C)$ by
\[
\psi_{1}(w)=\begin{cases}
(w^{1}-\sqrt{-1}w^{2}:w^{3}+\sqrt{-1}w^{4} ), & w=(w^{1}:w^{2}:w^{3}:w^{4}) \in \Q^{2}(\C)\backslash E \\
\lim_{u\not\in E, u\to w}\psi_{1}(w),  &  \mathrm{otherwise}\;   
\end{cases}
\]
where $E=\{(w^{1}:w^{2}:w^{3}:w^{4})\in \Q^{2}(\C)\: |\: w^{1}-\sqrt{-1}w^{2}=w^{3}+\sqrt{-1}w^{4}=0\}$. 
Since we can check
the value $\lim_{u\not\in E, u\to w}\psi_{1}(w)$ exists in $\mathbb{P}^{1}(\C)$, we define $\phi$. 
Similarly, for each $w=(w^{1}:w^{2}:w^{3}:w^{4})\in\Q^{2}(\C)$ we define 
\[
\psi_{2}(w)=\begin{cases}
(w^{1}-\sqrt{-1}w^{2}:-w^{3}+\sqrt{-1}w^{4} ), & w \in \Q^{2}(\C)\backslash E' \\
\lim_{u\not\in E', u\to w}\psi_{2}(w),  &  \mathrm{otherwise}\;   
\end{cases}
\]
where $E=\{(w^{1}:w^{2}:w^{3}:w^{4})\in \Q^{2}(\C)\: |\: w^{1}-\sqrt{-1}w^{2}=w^{3}-\sqrt{-1}w^{4}=0\}$. 
By using these maps, we define the map 
$\Psi=(\psi_{1}, \psi_{2})\colon \Q^{2}(\C)\to \mathbb{P}^{1}(\C)\times \mathbb{P}^{1}(\C)$. 
If we consider the map $\Psi^{\ast}\colon \mathbb{P}^{1}(\C)\times \mathbb{P}^{1}(\C)\to \Q^{2}(\C)$ defined by 
\[
\Psi^{\ast}((z:w), (u,v))=(zu+wv:\sqrt{-1}(zu-wv):wu-zv:-\sqrt{-1}(wu+zv)),
\]
we can easily check that $\Psi^{\ast}\circ\Psi$ and $\Psi\circ\Psi^{\ast}$ are both identity maps. 
Therefore, $\Psi$ is bijective and so the quadric $\Q^{2}(\C)$ is biholomorphic with $\mathbb{P}^{1}(\C)\times \mathbb{P}^{1}(\C)$.

Let $x\colon M\to \R^{4}$ be a surface. Take a complex local coordinate $z$ on $M$, 
we set $\phi_{i}=({\partial x^{i}}/{\partial z})dz$ $(i=1,\ldots,4)$ and define the map 
\[
g=(g_{1}, g_{2})=((\phi_{1}-\sqrt{-1}\phi_{2}:\phi_{3}+\sqrt{-1}\phi_{4}), (\phi_{1}-\sqrt{-1}\phi_{2}:-\phi_{3}+\sqrt{-1}\phi_{4}))\:.
\]
Instead of the Gauss map $g\colon M\to\Q^{2}(\C)$, we consider the map $g\colon M\to\mathbb{P}^{1}(\C)\times\mathbb{P}^{1}(\C)$,
which we call the Gauss map of $M$ in the following. 

Next, we explain the Enneper-Weierstrass representation theorem for minimal surfaces in $\R^{4}$. 
\begin{theorem}
Let $x=(x^{1}, x^{2}, x^{3}, x^{4})\colon M\to \R^{4}$ be a non-flat minimal surface immersed in $\R^{4}$. 
Consider the holomorphic $1$-forms $\phi_{1}$, $\phi_{2}$, $\phi_{3}$, $\phi_{4}$ which is defined by
\[
\phi_{i}=\frac{\partial x^{i}}{\partial z}dz \: (i=1,\ldots,4)\;.
\]
and the holomorphic $1$-form and the meromorphic functions which is defined by 
\begin{equation}\label{ewrep}
hdz=\phi_{1}-\sqrt{-1}\phi_{2},\quad 
g_{1}=\frac{\phi_{3}+\sqrt{-1}\phi_{4}}{\phi_{1}-\sqrt{-1}\phi_{2}},\quad g_{2}=\frac{-\phi_{3}+\sqrt{-1}\phi_{4}}{\phi_{1}-\sqrt{-1}\phi_{2}}
\end{equation}
Then,
\begin{enumerate}
\item [(i)] we have 
\begin{equation}\label{ewrep2}
\left\{
\begin{array}{l}
\phi_{1} = \frac{1}{2}(1+g_{1}g_{2})hdz \,\, , \\
\phi_{2} = \frac{\sqrt{-1}}{2} (1-g_{1}g_{2})hdz\,\, , \\
\phi_{3} = \frac{1}{2}(g_{1}-g_{2})hdz\,\, ,   \\
\phi_{4} = -\frac{\sqrt{-1}}{2}(g_{1}+g_{2})hdz\,\, 
\end{array}
\right.
\end{equation}
and we recover the immersion $x$ by the real Abel-Jacobi map
\begin{equation}\label{AbelJacobi}
x(z)=\Re\displaystyle\int^{z}_{z_{0}} (\phi_{1}, \phi_{2}, \phi_{3}, \phi_{4})
\end{equation}
up to translation.
\item [(ii)] the metric induced from the standard metric in $\R^{4}$ is given by
\begin{equation}\label{metric4}
ds^{2}=\frac{1}{4} |h|^{2}(1+|g_{1}|^{2})(1+|g_{2}|^{2})|dz|^{2}\,\, .
\end{equation}
\item [(iii)] the zeros of $hdz$ of order $k$ coincide exactly with the poles $g_{1}$ or $g_{2}$ of order $k$. 
(We call it ``the regularity condition''). 
\end{enumerate}
\end{theorem}
We also can show the following restatement of Theorem \ref{minimal3} for $n=4$. 
\begin{theorem}
Let $M$ be an open Riemann surface, $hdz$ a non-zero holomorphic $1$-form and $g_{1}$ and $g_{2}$ are meromorphic 
functions on $M$. Assume that the zeros of $hdz$ of order $k$ coincide exactly with the poles $g_{1}$ or $g_{2}$ of order $k$ 
and the holomorphic $1$-forms $\phi_{1}$, $\phi_{2}$, $\phi_{3}$, $\phi_{4}$ defined by (\ref{ewrep2}) have no real 
periods. Then, for the functions $x^{1}$, $x^{2}$, $x^{3}$, $x^{4}$ defined by (\ref{AbelJacobi}), the surface 
\[
x=(x^{1}, x^{2}, x^{3}, x^{4})\colon M\to\R^{4}
\]  
is a minimal surface immersed in $\R^{4}$ whose Gauss map is the map $g=(g_{1}, g_{2})$ and whose induced metric 
is given by (\ref{metric4}).
\end{theorem}

Now the Gauss curvature $K$ of $M$ is given by 
\[
K=-\frac{8}{|h|^{2}(1+|g_{1}|^{2})(1+|g_{2}|^{2})}\Biggl(\frac{|g'_{1}|^{2}}{(1+|g_{1}|^{2})^{2}}+\frac{|g'_{2}|^{2}}{(1+|g_{2}|^{2})^{2}}\Biggr)
\] 
and the total curvature by 
\[
\tau(M)=\int_{M}KdA=-\int_{M}\Biggl(\frac{2|g_{1}'|^{2}}{(1+|g_{1}|^{2})^{2}}+\frac{2|g_{2}'|^{2}}{(1+|g_{2}|^{2})^{2}}\Biggr)|dz|^{2}
\]
where $dA$ is the area form of $M$. When the total curvature of a complete minimal surface is finite, 
the surface is called \textit{an algebraic minimal surface}. The following theorem is the restatement 
for Theorem \ref{Chern-Osserman} for $n=4$.
\begin{theorem}[Huber-Osserman]\label{ho}
An algebraic minimal surface $x\colon M\to \R^{4}$ satisfies the followings :
\begin{enumerate}
\item [(i)] $M$ is conformally equivalent to $\overline{M}\backslash \{p_{1},\ldots,p_{k}\}$ where $\overline{M}$ 
is a compact Riemann surface, and $p_{1},\ldots,p_{k}$ are finitely many points of $\overline{M}$.
\item [(ii)] The Weierstrass data $(hdz, g_{1}, g_{2})$ extend meromorphically to $\overline{M}$.
\end{enumerate}
\end{theorem}
For a complete minimal surface in $\R^{4}$, the definition of ``pseudo-algebraic'' is as follows.

\begin{definition}\label{pa4}
We call a complete minimal surface in $\R^{4}$ \textit{pseudo-algebraic}, if the following conditions are satisfied:
\begin{enumerate} 
\item [(i)] The Weierstrass data $(hdz, g_{1}, g_{2})$ is defined on a Riemann surface $M=\overline{M}\backslash \{p_{1},\ldots,p_{k}\}$ where $\overline{M}$ 
is a compact Riemann surface, and $p_{1},\ldots,p_{k}\in\overline{M}$.
\item [(ii)] The Weierstrass data $(hdz, g_{1}, g_{2})$ can extend meromorphically to $\overline{M}$.
\end{enumerate}
We call $M$ \textit{the basic domain} of the pseudo-algebraic minimal surface under consideration. 
\end{definition}

Algebraic minimal surfaces are certainly pseudo-algebraic. 
The following examples are also pseudo-algebraic.

\begin{example}[Mo-Osserman \cite{MO}]\label{e1}
Let $M=\C\backslash\{a_{1},a_{2},a_{3}\}$ for distinct $a_{1}, a_{2}, a_{3}\in \C$, 
the Weierstrass data is defined on $M$ 
by
\[
(hdz,g_{1},g_{2})=\Bigl(\frac{dz}{\prod_{i=1}^{3}(z-a_{i})},z,z\Bigr)\,\, .
\]
As this data does not satisfy period condition, we get a minimal surface $x\colon\D\to \R^{4}$ 
on the universal covering 
disk of $M$. In particular, it has infinite total curvature. 
We can see that the surface is complete and both $g_{1}$ and $g_{2}$ 
omit four values $a_{1}, a_{2}, a_{3}, \infty$. This surface does not lie fully in $\R^{4}$ 
because the component function 
$x^{3}$ is equal to $0$. Thus it is one-degenerate. (For details, see \cite{HO1})
\end{example}  
\begin{example}[Mo-Osserman \cite{MO}]\label{e2}
Let $M=\C\backslash\{a_{1},a_{2}\}$ for distinct $a_{1}, a_{2} \in \C$, 
the Weierstrass data is defined on $M$ 
by
\[
(hdz,g_{1},g_{2})=\Bigl(\frac{dz}{\prod_{i=1}^{2}(z-a_{i})},z,0\Bigr)\,\, .
\]
As this data does not satisfy period condition, we also get a complete minimal surface with infinite total curvature 
$x\colon\D\to \R^{4}$ on the universal covering disk of $M$. 
We can see that $g_{1}$ omits three values 
$a_{1}, a_{2}, \infty$. This surface is a complex curve in $\C^{2}\simeq \R^{4}$ 
because the Gauss map $g_{2}$ is constant.
\end{example}  

\subsection{Ramification estimate and unicity theorem}
We give bound estimates for the totally ramified value number of the Gauss map of pseudo-algebraic minimal 
surfaces in $\R^{4}$. 
\begin{theorem}[Kawakami, \cite{Ka2}]\label{main1}
Consider a non-flat pseudo-algebraic minimal surface in $\R^{4}$ with the basic domain $M=\overline{M}\backslash \{p_{1},\ldots,p_{k}\}$. 
Let $G$ be the genus of $\overline{M}$, $d_{i}$ be the degree of $g_{i}$ considered as a map $\overline{M}$ 
and $\nu_{g_{i}}$ be the totally ramified value number of $g_{i}$. 
\begin{enumerate}
\item [(i)]In the case $g_{1}\not\equiv const.$ and $g_{2}\not\equiv const.$, then $\nu_{g_{1}}\leq 2$, or $\nu_{g_{2}}\leq 2$, or
\begin{equation}\label{eq1}
\frac{1}{\nu_{g_{1}}-2}+\frac{1}{\nu_{g_{2}}-2}\geq R_{1}+R_{2}\geq 1, \quad R_{i}=\frac{d_{i}}{2G-2+k}\:(i=1, 2) 
\end{equation}
and for an algebraic minimal surface, $R_{1}+R_{2}> 1$\,\, .
\item [(ii)]In the case where one of $g_{1}$ and $g_{2}$ is costant, say $g_{2}\equiv const.$, then 
\begin{equation}\label{eq2}
\nu_{g_{1}}\le 2+\frac{1}{R_{1}}, \quad R_{1}=\frac{d_{1}}{2G-2+k}\geq 1
\end{equation}
and for an algebraic minimal surface, $R_{1}>1$\,\, .
\end{enumerate}
\end{theorem} 
\begin{proof}
By a suitable rotation of the surface, we may assume that both $g_{1}$ and $g_{2}$ are no pole at $p_{j}$, 
and have only simple poles. By the completeness, $hdz$ has poles of order $\mu_{j}\geq 1$ at 
$p_{j}$. By Proposition \ref{order2}, the period condition implies $\mu_{j}\geq 2$, however here we do not assume this. Let $\alpha_{s}$ be (simple) poles 
of $g_{1}$, $\beta_{t}$ (simple) poles of $g_{2}$. The following table shows the relation between zeros and poles of $g_{1}$, 
$g_{2}$ and $hdz$. The upper index means the order.
\begin{center}
\begin{tabular}{|c|c|c|c|}\hline
$z$ & $\alpha_{s}$  & $\beta_{t}$ & $p_j$  \\\hline
$g_{1}$ & $\infty^{1}$ &   &   \\\hline
$g_{2}$ &   & $\infty^{1}$ &   \\\hline
$hdz$   & $0^1$ & $0^1$ & $\infty^{\mu_i}$ \\\hline
\end{tabular}
\end{center}
Applying the Riemann-Roch formula to the meromorphic differential $hdz$ on $\overline{M}$, we obtain
\[
d_{1}+d_{2}-\displaystyle \sum_{i=1}^{k}\mu_{i}=2G-2\,\, .
\]
Note that this equality depends on above setting of poles of $g_{1}$ and $g_{2}$, though $d_{1}$ and $d_{2}$ are 
invariant. Thus we get
\begin{equation}\label{eq3}
d_{1}+d_{2}=2G-2+\displaystyle \sum_{i=1}^{k}\mu_{i}\ge 2G-2+k
\end{equation}
and
\begin{equation}\label{eq4}
R_{1}+R_{2}=\frac{d_{1}+d_{2}}{2G-2+k}\ge 1\,\, .
\end{equation}
When $M$ is an algebraic minimal surface, we have $\mu_{j}\geq 2$ and so $R_{1}+R_{2}>1$.

Now we prove (i). Assume $g_{i}$ is not constant and omits $r_{i0}$ values. Let $n_{i0}$ be the sum of the branching 
orders of $g_{i}$ at the inverse image of these exceptional values. We see
\begin{equation}\label{eq5}
k\ge d_{i}r_{i0}-n_{i0}.
\end{equation}
Let $b_{i1},\ldots,b_{il_{0}}$ be the totally ramified values which are not exceptional values, and $n_{ir}$ the 
sum of branching orders of $g_{i}$ at the inverse image of these values. For each $b_{ij}$, we denote 
\[
\nu_{ij}=min_{g^{-1}(b_{ij})}\{\text{multiplicity of }g(z)=b_{ij} \}\,\, ,
\]
then the number of points in the inverse image 
$g_{i}^{-1}(b_{ij})$ is less than or equal to $d_{i}/\nu_{ij}$. Thus we obtain
\begin{equation}\label{eq6}
d_{i}l_{0}-n_{ir}\le \displaystyle\sum_{j=1}^{l_{0}}\frac{d_{i}}{\nu_{ij}}\,\, .
\end{equation}
This implies
\begin{equation}\label{eq7}
l_{0}-\sum_{j=1}^{l_{0}}\frac{1}{\nu_{ij}} \leq \frac{n_{ir}}{d_{i}}\,\, .
\end{equation}
Let $n_{i1}$ be the total branching order of $g_{i}$. Then applying Riemann-Hurwitz's theorem to the meromorphic function 
$g_{i}$ on $\overline{M}$, we obtain
\begin{equation}\label{eq8}
n_{i1}=2(d_{i}+G-1)\,\, .
\end{equation}
By (\ref{eq5}), (\ref{eq7}) and (\ref{eq8}), we get
\begin{equation}\label{eq9}
\nu_{g_{i}}=r_{i0}+\displaystyle\sum_{j=1}^{l_{0}}(1-\frac{1}{\nu_{ij}})\le \frac{n_{i0}+ k}{d_{i}}+\frac{n_{ir}}{d_{i}}\le \frac{n_{i1}+k}{d_{i}}=2+\frac{1}{R_{i}}\,\, .
\end{equation}
When $\nu_{g_{1}}>2$ and $\nu_{g_{2}}>2$, 
\[
\frac{1}{\nu_{g_{i}}-2}\ge R_{i}\quad(i=1,2)\,\, .
\]
Hence we get
\[
\frac{1}{\nu_{g_{1}}-2}+\frac{1}{\nu_{g_{2}}-2}\ge R_{1}+R_{2}\,\, .
\]

Next, we show (ii). Then, the simple poles of $g_{1}$ coincides exactly with the simple zeros of $hdz$ 
and $hdz$ has a pole of order $\mu_{j}$ at each $p_{j}$. Applying the Riemann-Roch formula to the meromorphic differential $hdz$ on $\overline{M}$, we obtain
\[
d_{1}-\displaystyle \sum_{i=1}^{k}\mu_{i}=2G-2\,\, .
\]
Thus we get
\begin{equation}\label{eq10}
d_{1}=2G-2+\displaystyle \sum_{i=1}^{k}\mu_{i}\ge 2G-2+k
\end{equation}
and
\begin{equation}\label{eq11}
R_{1}=\frac{d_{1}}{2G-2+k}\ge 1\,\, .
\end{equation}
When $M$ is an algebraic minimal surface, we have $\mu_{j}\geq 2$ and so $R_{1}>1$. 
By (\ref{eq9}), we get
\[
\nu_{g_{1}} \leq 2+\frac{1}{R_{1}}\,\, .
\]
Thus, we complete the proof of this theorem.
\end{proof}

We have the following result as an immediate consequence of Theorem \ref{main1}.
\begin{corollary}[cf.\:Fujimoto\cite{F1}, Hoffman-Osserman\cite{HO1}] \label{co1}
Let $x\colon M\to\R^{4}$ be a pseudo-algebraic minimal surface, $g=(g_{1}, g_{2})$ be its Gauss map.
\begin{enumerate}
\item [(i)]In the case $g_{1}\not\equiv const.$ and $g_{2}\not\equiv const.$, if both $g_{1}$ and $g_{2}$ omit more than four values, 
then $M$ must be a plane. In particular, if $M$ is an algebraic minimal surface and if both $g_{1}$ and $g_{2}$ omit more than three values, 
then $M$ must be a plane.
\item [(ii)]In the case where one of $g_{1}$ and $g_{2}$ is constant, say $g_{2}\equiv const.$, if $g_{1}$ omits more than three values, 
then $M$ must be a plane. In particular, if $M$ is an algebraic minimal surface and if $g_{1}$ omits more than two values, then 
$M$ must be a plane.  
\end{enumerate}
\end{corollary}
Example \ref{e1} and Example \ref{e2} show Corollary \ref{co1} is the best possible for pseudo-algebraic minimal surfaces. 
The following example shows Corollary \ref{co1} (ii) is the best possible also for algebraic minimal surfaces.
\begin{example}[Kawakami, \cite{Ka2}]\label{e3}
Let $M=\C\backslash\{0\}$, the Weierstrass data is defined on $M$ 
by
\[
(hdz,g_{1},g_{2})=\Bigl(\frac{dz}{z^{3}},z,c\Bigr)\,\,.
\]
where $c$ is a complex number. As this data satisfy the regularity condition, the period condition on $M$ and the surface is complete, 
we get an algebraic minimal surface $x\colon M\to \R^{4}$ whose Gauss map $g_{1}$ omits two values $0, \infty$.
\end{example}
\begin{remark}
In Section 5, We state the results on remification estimate for the Gauss map $g\colon M \to \mathbb{P}^{n-1}(\C)$ of 
a complete minimal surface in $\R^{n}$. However these results do not cover Theorem \ref{main1} and Corollary \ref{co1}
because corresponding hyperplanes in $\mathbb{P}^{3}(\C)$ are not necessary located in general position. 
(For details, see [14, p353].)
\end{remark} 

Next, we give the unicity theorem for the Gauss map of a pseudo-algebraic minimal surface in $\R^{4}$.
\begin{theorem}[Kawakami, \cite{Ka2}]\label{main2}
Consider two non-flat pseudo-algebraic minimal surfaces in $\R^{4}$, $M_{A}$ and $M_{B}$ 
with the same basic domain $M=\overline{M}\backslash \{p_{1},\ldots,p_{k}\}$. 
Let $G$ be the genus of $\overline{M}$, and $g_{A}=(g_{A1}, g_{A2})$, $g_{B}=(g_{B1}, g_{B2})$ be the Gauss map 
of $M_{A}$ and $M_{B}$ respectively. For each $i$ $(i=1, 2)$, assume that both $g_{Ai}$ and $g_{Bi}$ have the same degree $d_{i}$ when considered as a map on $\overline{M}$.
\begin{enumerate} 
\item [(i)] In the case $g_{A1}\not\equiv g_{B1}$ and $g_{A2}\not\equiv g_{B2}$, let $a_{1},\ldots, a_{p}\in \hat{C}$, 
$b_{1},\ldots, b_{q}\in \hat{\C}$ be distinct points such that $g_{A1}^{-1}(a_{j})\cap M = g_{B1}^{-1}(a_{j})\cap M$ 
for $1\leq j\leq p$, $g_{A2}^{-1}(b_{k})\cap M = g_{B2}^{-1}(b_{k})\cap M$ for $1\leq k\leq q$ respectively. 
If $p > 4$ and $q > 4$, then
\begin{equation}\label{eq12}
\frac{1}{p-4}+\frac{1}{q-4}\ge R_{1}+R_{2} \geq 1,\quad R_{i}=\frac{d_{i}}{2G-2+k}\ (i=1,2)\,\, .
\end{equation}
In particular, if $p\geq 7$ and $q\geq 7$ then $g_{A}\equiv g_{B}$. 
\item [(ii)] In the case $g_{A1}\not\equiv g_{B1}$ and $g_{A2}\equiv g_{B2}\equiv const.$, 
let $a_{1},\ldots, a_{p}\in \hat{\textbf{C}}$ 
be distinct points such that $g_{A1}^{-1}(a_{j})\cap M = g_{B1}^{-1}(a_{j})\cap M$ for $1\leq j\leq p$. Then
\begin{equation}\label{eq13}
p\le 4+\frac{1}{R_{1}},\quad R_{1}=\frac{d_{1}}{2G-2+k}\,\, .
\end{equation}
In particular, if $p\geq 6$ then $g_{A}\equiv g_{B}$.
\end{enumerate}
\end{theorem}
\begin{proof}
Put
\[
\delta_{j}=\sharp(g_{A1}^{-1}(a_{j})\cap M)=\sharp(g_{B1}^{-1}(a_{j})\cap M)
\]
where $\sharp$ denotes the number of points. Then we have 
\begin{equation}\label{eq14}
pd_{1}\le \displaystyle\sum_{j=1}^{p}\delta_{j}+n_{11}+k
\end{equation}
using the same notation as in proof of Theorem \ref{main1}. 
Consider a meromorphic function $\varphi= 1/(g_{A1}-g_{B1})$ on $M$. $\varphi$ has a pole, 
while the total number of the poles of $\varphi$ is at most $2d_{1}$, we get
\begin{equation}\label{eq15}
\displaystyle\sum_{j=1}^{p}\delta_{j}\le 2d_{1}\,\, .
\end{equation}
Then from (\ref{eq14}) and (\ref{eq15}), we obtain 
\[
pd_{1}\le 2d_{1}+n_{11}+k
\]
and
\begin{equation}\label{eq16}
p\le\frac{2d_{1}+n_{11}+k}{d_{1}}=4+\frac{1}{R_{1}}\,\, .
\end{equation}
Similarly we obtain 
\begin{equation}\label{eq17}
q\leq 4+\frac{1}{R_{2}}\,\, .
\end{equation}
From (\ref{eq16}) and (\ref{eq17}), we get (\ref{eq12}) and (\ref{eq13}) immediately.
\end{proof}

We give an example which shows $(p, q)=(7, 7)$ in Theorem \ref{main2} (i) is the best possible
for pseudo-algebraic minimal surfaces.
\begin{example}[Kawakami, \cite{Ka2}]
Taking a complex number $\alpha$ with $\alpha\not=0, \pm 1$, we consider the Weierstrass data 
\[
(hdz,g_{1},g_{2})=\Bigl(\frac{dz}{z(z-\alpha)(\alpha z-1)},z,z\Bigr)
\]
and the universal covering surface $M$ of $\C\backslash \{0, \alpha, 1/\alpha\}$. Then we can constract
a pseudo-algebraic minimal surface on $M$. On the other hand, 
we can constract a pseudo-algebraic minimal surface on $M$ whose Weierstrass data is 
\[
(hdz,\bar{g_{1}},\bar{g_{2}})=\Bigl(\frac{dz}{z(z-\alpha)(\alpha z-1)},\frac{1}{z},\frac{1}{z}\Bigr)\,\, .
\]
For the maps $g_{i}$ and $\bar{g_{i}}$, we have $g_{i}\not\equiv \bar{g_{i}}$ and $g_{i}^{-1}(\alpha_{j})=
\bar{g_{i}}^{-1}(\alpha_{j})$ for six values
\[
\alpha_{1}:=0,\:\alpha_{2}:=\infty,\:\alpha_{3}:=\alpha,\:\alpha_{4}:=\frac{1}{\alpha},\:\alpha_{5}:=1,\:
\alpha_{6}:=-1.
\]
\end{example}
We also give an example which shows $p=6$ in Theorem \ref{main2} (ii) is the best possible
for pseudo-algebraic minimal surfaces.
\begin{example}[Kawakami, \cite{Ka2}]
Taking a complex number $\alpha$ with $\alpha\not=0, \pm 1$, we consider the Weierstrass data 
\[
(hdz,g_{1},g_{2})=\Bigl(\frac{dz}{z(z-\alpha)},z,0\Bigr)
\]
and the universal covering surface $M$ of $\C\backslash \{0, \alpha\}$. Then we can constract
a pseudo-algebraic minimal surface on $M$. On the other hand, 
we can constract a pseudo-algebraic minimal surface on $M$ whose Weierstrass data is 
\[
(hdz,\bar{g_{1}},\bar{g_{2}})=\Bigl(\frac{dz}{z(z-\alpha)},\frac{1}{z},0\Bigr)\,\, .
\]
For the maps $g_{1}$ and $\bar{g_{1}}$, we have $g_{1}\not\equiv \bar{g_{1}}$ and $g_{1}^{-1}(\alpha_{j})=
\bar{g_{1}}^{-1}(\alpha_{j})$ for five values
\[
\alpha_{1}:=0,\:\alpha_{2}:=\infty,\:\alpha_{3}:=\alpha,\:\alpha_{4}:=1,\:\alpha_{5}:=-1.
\]
\end{example}

\section{The Gauss map of pseudo-algebraic minimal surfaces in $\R^{n}$}
In this section, we shall study the 
Gauss map of pseudo-algebraic minimal surfaces in $\R^{m}$.

\subsection{Some results of algebraic curves in the projective space}
In this subsection, we give some resuls on a holomorphic map of a compact Rimann surface with genus $G$ 
(it is denoted by $\overline{M}_{G}$)
into $\mathbb{P}^{n}(\C)$ (we call it ``{\em algebraic curve}'') to show ramification estimates for the Gauss map. 

First, we recall some results on algebraic curve. 
\begin{definition}
An algebraic curve $f\colon \overline{M}_{G}\to \mathbb{P}^{n}(\C)$ is said to be {\em linearly nondegenerate} if the image of $f$ is 
not included in any hyperplane in $\mathbb{P}^{n}(\C)$. 
\end{definition}

Assume that $f\colon \overline{M}_{G}\to \mathbb{P}^{n}(\C)$ is an algebraic curve. For a fixed homogeneous coordinates 
$(w^{0}:\cdots:w^{n})$ we set 
\[
V_{i}=\{(w^{0}:\cdots:w^{n})\, |\, w^{i}\not=0\}\quad (0\leq i\leq n)\;.
\]
Then, every $p\in\overline{M}_{G}$ has a neighborhood $U$ of $p$ such that $f(U)\subset V_{i}$ for some $i$ and $f$ 
has a representation 
\[
f=(f_{0}:\cdots:f_{i-1}:1:f_{i+1}:\cdots:f_{n})
\]
on $U$ with holomorphic functions $f_{0},\ldots,f_{i-1},f_{i+1},\ldots,f_{n}$.

\begin{definition}
For an open subset $U$ of $\overline{M}_{G}$ we call a representation $f=(f_{0}:\cdots:f_{n})$ to be a {\em reduced representation} 
of $f$ on $U$ if $f_{0},\ldots,f_{n}$ are holomorphic functions on $U$ and have no common zero. 
\end{definition}

Let $f\colon \overline{M}_{G}\to\mathbb{P}^{n}(\C)$ be a linearly nondegerate algebraic curve. Take a point $p\in \overline{M}$.
For a suitable choice of homogeneous coordinates $(w^{0}:\cdots:w^{n})\in \mathbb{P}^{n}(\C)$, the equation of 
the curve can be put locally into the normal form
\begin{equation}\label{deldel}
(w^{0}:\cdots:w^{n})=(z^{\delta_{0}}+\cdots: \cdots :z^{\delta_{n}}+\cdots)\;, 
\end{equation}
where
\[
0=\delta_{0}<\delta_{1}<\cdots<\delta_{n}
\]
and $z$ is a complex local coordinate with $z(p)=0$ on $\overline{M}_{G}$. The integers 
\[
\nu_{i}={\delta}_{i+1}-{\delta}_{i}-1,\; 0\leq i\leq n-1
\]
are called the {\em stationary induces of order i} at the point $z=0$ (Geometrically, 
this is the order of the associated curve of rank $i$, i.e., the curve formed by the osculating spaces of 
dimension $i$). And we have 
\begin{equation}\label{stationary}
\displaystyle \sum_{0\leq i\leq n-1} (n-i)\nu_{i}(p)+\frac{1}{2}n(n+1)={\delta}_{1}(p)+\cdots+\delta_{n}(p)\;.
\end{equation}
The stationary points, i.e., points with a non-zero stationary index, are isolated and hence are finite in number. 
We will denote by $\sigma_{i}$ the sum of all stationary indices of order $i$. Then Pl\"ucker formula are 
\begin{equation}\label{Plucker}
\nu_{i-1}-2\nu_{i}+\nu_{i+1}=2(G-1)-\sigma_{i},\; 0\leq i\leq n-1,
\end{equation}
with the convention $\nu_{-1}=\nu_{n}=0$\:. Pl\"ucker formula is a generalization of Riemann-Hurwitz's theorem. 
For the proof of Pl\"ucker formula, see [\cite{F4}, p177]. From (\ref{Plucker}) it follows that 
\begin{equation}\label{Plucker2}
\displaystyle \sum_{0\leq i\leq n-1} (n-i)\sigma_{i}=(n+1)\deg(f)+n(n+1)\; ,
\end{equation}   
where $\deg(f)$ is the degree of $f$.
\begin{definition}
Let $H_{1},\ldots,H_{q}$ be hyperplanes in $\mathbb{P}^{n}(\C)$ and $L_{1},\ldots,L_{q}$ be the corresponding linear 
forms. We say that $H_{1},\ldots,H_{q}$ are {\em general position} if for any injective map $\mu\colon \{0, 1, \ldots, n\}
\to \{1,\ldots,q\}$, $L_{\mu(0)},\ldots,L_{\mu_{n}}$ are linearly independent.   
\end{definition}

\begin{theorem}[Chern-Osserman \cite{CO}, Jin and Ru \cite{JR}]\label{SMT1}
Let $\overline{M}_{G}$ be a compact Riemann surface with genus $G$ and let be $f\colon \overline{M}_{G}\to \mathbb{P}^{n}(\C)$ 
be a linearly nondegenerate algebraic curve. Let $H_{1},\ldots,H_{q}$ be the hyperplanes in $\mathbb{P}^{n}(\C)$, 
located in general position. Let $E=\cup_{j=1}^{q} f^{-1}(H_{j}) $. Then, 
\[
\{q-(n+1)\}\deg(f)\leq \frac{1}{2}n(n+1)\{2(g+1)+\sharp E\}\;,
\]
where $\sharp$ denotes the number of points.
\end{theorem}
\begin{proof}
We denote $E=\{p_{1},\ldots,p_{s}\}$. First of all, if $f(\overline{M}_{G})$ intersects $H_{j}$ at certain point 
$p_{l}\in E$ with some multiplicity $v_{p_{l}}(L_{j}(f))$, where $L_{j}$ is the linear form corresponding to $H_{j}$. 
Then, by the definition, for every $1\leq j\leq q$,
\begin{equation}\label{Plucker3}
\displaystyle \sum_{1\leq l\leq s} v_{p_{l}}(L_{j}(f))= \deg(f)\;.
\end{equation}
Secondly, since the hyperplanes $H_{1},\ldots,H_{q}$ are in general position, at most $n$ hyperplanes can 
intersect $f(\overline{M}_{G})$ at $p_{l}$, hence there exists subset $A\subset \{1, 2, \ldots, q\}$ with $\sharp A= n$ 
such that 
\begin{equation}\label{Plucker4}
\sum_{1\leq j\leq q} v_{p_{l}}(L_{j}(f))\leq \displaystyle \sum_{i\in A}v_{p_{l}}(L_{i}(f))\;.
\end{equation} 
Take a complex local coordinate $z$ for $\overline{M}_{G}$ at $p_{l}$ such that $z(p_{l})=0$. 
At $p_{l}$, the maximum possible value of $v_{p_{l}}(L_{j}(f))$, $i\in A$, is $\delta_{n}(p_{l})$, and for 
the unique hyperplane $w^{n}=0$\;. A second hyperplane can intersect $f(\overline{M}_{G})$ at $p_{l}$ with 
multiplicities at most $\delta_{n-1}(p_{l})$, etc. It follows that 
\begin{equation}\label{Plucker5}
\displaystyle \sum_{i\in A}v_{p_{l}}(L_{i}(f))\leq {\delta}_{1}(p_{l})+\cdots+{\delta}_{n}(p_{l})\;.
\end{equation}
By (\ref{stationary}), we get 
\[
{\delta}_{1}(p_{l})+\cdots+{\delta}_{n}(p_{l})=\displaystyle \sum_{0\leq i\leq n-1} (n-i)\nu_{i}(p_{l})+\frac{1}{2}n(n+1)\;.
\]
Combining this with (\ref{Plucker2}), (\ref{Plucker3}), (\ref{Plucker4}), and (\ref{Plucker5}), we get
\[
q\,\deg(f)\leq (n+1)\deg(f) +n(n+1)(G-1)+\frac{1}{2}n(n+1)\sharp E\;.
\]
\end{proof}

Now, we extend Theorem \ref{SMT1} to the degenerate case. Assume that $f\colon \overline{M}_{G}\to \mathbb{P}^{n}(\C)$ 
be an algebraic curve (not necessarily linearly nondegenerate) and $f(\overline{M}_{G})$ is contained in some $r$-dimensional 
projective subspace of $\mathbb{P}^{n}(\C)$, however it is not in any subspace of dimensional lower than $r$, where $1\leq r\leq n$. 
Then $f\colon \overline{M}_{G}\to \mathbb{P}^{r}(\C)$ is a linearly nondegenerate algebraic curve. 
Let $H_{1}, \ldots, H_{q}$ be the hyperplanes in $\mathbb{P}^{n}(\C)$, located in general position. 
Then their restrictions to $\mathbb{P}^{r}(\C)$, $H_{1}\cap \mathbb{P}^{r}(\C), \ldots, H_{q}\cap \mathbb{P}^{r}(\C)$ 
are in $n$-subgeneral position in $\mathbb{P}^{r}(\C)$, i.e. $n+1$ of them (regared as linear forms) span the 
($r+1$)-demensional complex plane $\C^{r+1}$. The difficulty of degenerate case is that hyperplanes $H_{1},\ldots,H_{q}$ in 
$\mathbb{P}^{n}(\C)$ in general position may not necessarily in general position after being restricted to $\mathbb{P}^{r}(\C)$. 
So we have to use the following techniques of Nochka to overcome this difficulty. 
These techniques are essential in the solution of the Cartan conjecture (cf.\:\cite{F4}) . 
\begin{theorem}[Nochka, \cite{F4}, \cite{R3}]\label{Nochka1}
Let $H_{1}, \ldots, H_{q}$ be hyperplanes in $\mathbb{P}^{r}(\C)$ in $n$-subgeneral position with $2n-r+1\leq q$. 
Let $L_{1}, \ldots, L_{q}$ be the corresponding linear forms. Then there exists a function $\omega\colon \{1, \ldots ,q\}\to (0, 1]$ 
called a Nochka weight and a real number $\theta\geq 1$ called Nochka constant satisfying the following properties: 
\begin{enumerate}
\item [(i)] If $j\in \{1, \ldots, q\}$, then $0\leq \omega(j)\theta\leq 1$. 
\item [(ii)] $q-2n+r-1=\theta(\sum_{j=1}^{q}\omega(j)-r-1)$\;.
\item [(iii)]If $\emptyset \not= B \subset \{1, \ldots, q\}$ with $\sharp B \leq n+1$, then $\sum_{j\in B} \omega(j)\leq \dim L(B)$, 
where $L(B)$ is the linear space generated by $\{L_{j}\, |\, j\in B\}$. 
\item [(iv)] $1\leq (n+1)/(r+1)\leq \theta \leq (2n-r+1)/(r+1)$\;.
\item [(v)] Given real numbers $\lambda_{1},\ldots, \lambda_{q}$ with $\lambda_{j}\geq 1$ for $1\leq j\leq q$, 
and given any $Y\subset \{1,\ldots ,q\}$ with $0<\sharp Y\leq n+1$, there exists a subset $M$ of $Y$ with $\sharp M=\mathrm{dim}\:L(Y)$ such that 
$\{L_{j}\}_{j\in M}$ is a basis for $L(Y)$ where $L(Y)$ is the linear space generated by $\{L_{j}\:|\: j\in Y\}$, and 
\[
\displaystyle \prod_{j\in Y}\lambda_{j}^{\omega(j)}\leq \displaystyle \prod_{j\in M}\lambda_{j}\;.
\]
\end{enumerate}
\end{theorem}

\begin{theorem}[Jin and Ru \cite{JR}]\label{SMT2}
Let $f\colon \overline{M}_{G}\to \mathbb{P}^{n}(\C)$ be a non-constant algebraic curve. Assume that $f(\overline{M}_{G})$ 
is contained in some $r$-dimensional projective subspace of $\mathbb{P}^{n}(\C)$, however it is not in any subspace of 
dimension lower than $r$, where $1\leq r\leq n$. Let $H_{1}, \ldots, H_{q}$ be the hyperplanes in $\mathbb{P}^{n}(\C)$, 
located in general position. Let $E=\cup_{j=1}^{q} f^{-1}(H_{j}) $. Then 
\[
(q-2n+r-1)\deg(f)\leq \frac{1}{2}r(2n-r+1)\{2(G-1)+\sharp E\}\;.
\]
\end{theorem}
\begin{proof}
By the assumption, $f\colon \overline{M}_{G}\to \mathbb{P}^{k}(\C)$ is linearly non-degenerate. Since $H_{1}, \ldots, H_{q}$ are 
general position, their restrictions to $\mathbb{P}^{r}(\C)$, $H_{1}\cap \mathbb{P}^{r}(\C), \ldots, H_{q}\cap \mathbb{P}^{r}(\C)$ 
are in $n$-subgeneral position in $\mathbb{P}^{r}(\C)$. For simplicity, we still denote $H_{j}\cap \mathbb{P}^{r}(\C)$ 
as $H_{j}$, $1\leq j\leq q$. Let $L_{j}$ be the linear forms defining $H_{j}$, $1\leq j\leq q$. Let $\omega(j)$ be the 
Nochka weights associated to the hyperplanes $H_{j}$, $1\leq j\leq q$. We denotes $E=\{p_{1},\ldots,p_{s}\}$. 
For any point $p_{l}\in E$, taking ${\lambda}_{j}=e^{{v}_{p_{l}}({L}_{j}(f))}$, and using (v) of Theorem \ref{Nochka1}, 
there exists $L_{{P}_{l}, 1}, \ldots, L_{{P}_{l}, r}$ such that they are linearly independent, and that 
\[
\displaystyle \prod_{j=1}^{q} e^{\omega(j){v}_{p_{l}}({L}_{j}(f))} \leq \displaystyle \prod_{j=1}^{r} e^{{v}_{p_{l}}({L}_{p_{l}, j}(f))}\;. 
\] 
This gives 
\begin{equation}\label{Nochka2}
\displaystyle \sum_{j=1}^{q} \omega(j) {v}_{p_{l}}({L}_{j}(f))\leq \displaystyle \sum_{j=1}^{r} {v}_{p_{l}}({L}_{p_{l}, j}(f))\;. 
\end{equation}
In the same way as (\ref{Plucker5}), we have
\begin{equation}\label{Nochka2}
\displaystyle \sum_{j=1}^{r} {v}_{p_{l}}({L}_{p_{l}, j}(f))\leq {\delta}_{1}(p_{l})+\cdots+{\delta}_{n}(p_{l}) \;. 
\end{equation}
By (\ref{stationary}), we get 
\begin{equation}\label{Nochka3}
{\delta}_{1}(p_{l})+\cdots+{\delta}_{n}(p_{l}) = \displaystyle \sum_{0\leq i\leq r-1} (r-i){\nu}_{i}({P}_{l})+\frac{1}{2}r(r+1) \;.
\end{equation}
Thus we have 
\begin{equation}\label{Nochka4}
\displaystyle \sum_{j=1}^{q} \displaystyle \sum_{l=1}^{s} \omega(j) {v}_{p_{l}}({L}_{j}(f))\leq (r+1)\deg(f)+\frac{r(r+1)}{2}(2G-2+\sharp E)\;.
\end{equation}
By (\ref{Plucker2}), we have 
\begin{equation}
\displaystyle \sum_{j=1}^{q} \omega(j)\deg(f)\leq (r+1)\deg(f)+\frac{r(r+1)}{2}(2G-2+\sharp E)\;.
\end{equation}
Using (ii) and (iv) of Theorem \ref{Nochka1}, we get
\[
(q-2n+r-1)\deg(f)\leq \frac{\theta k(k+1)}{2}(2G-2+\sharp E)\leq \frac{r(2n-r+1)}{2}(2G-2+\sharp E)\;.
\]
We have thus proved the theorem. 
\end{proof}

In the following theorem, we modify Theorem $\ref{SMT2}$ to the case that $E$ is an arbitrary finite subset of $\overline{M}$. 

\begin{theorem}[Jin and Ru \cite{JR}]\label{SMT3}
Let $f\colon \overline{M}_{G} \to \mathbb{P}^{n}(\C)$ be a non-constant algebraic curve. 
Assume that $f(\overline{M}_{G})$ 
is contained in some $r$-dimensional projective subspace of $\mathbb{P}^{n}(\C)$, however it is not in any subspace of 
dimension lower than $r$, where $1\leq r\leq n$. Let $H_{1}, \ldots, H_{q}$ be the hyperplanes in $\mathbb{P}^{n}(\C)$, 
located in general position and let $L_{1}, \ldots, L_{q}$ be the corresponding linear forms. 
Let $E$ be a finite subset of $\overline{M}_{G}$. Then
\[
(q-2n+r-1)\deg(f)\leq \displaystyle \sum_{j=1}^{q} \displaystyle \sum_{p \not\in E}\min\{r, v_{p}(L_{j}(f))\} 
+\frac{1}{2}r(2n-r+1)\{2(G-1)+\sharp E\}\:,
\]
where $v_{p}(L_{j}(f))$ is the vanishing order of $L_{j}(f)$ at the point $p$. 
\end{theorem}
\begin{proof}
The above inequality trivially holds for $q \leq 2n-r+1$\;. So we assume that $q > 2n-r+1$\;. 
By the assumeption, $f\colon \overline{M}_{G}\to \mathbb{P}^{r}(\C)$ is  linearly nondegenerate. 
Since $H_{1}, \ldots, H_{q}$ are in general position, their restriction (to $\mathbb{P}^{r}(\C)$) 
$H_{1}\cap \mathbb{P}^{r}(\C), \ldots, H_{q}\cap \mathbb{P}^{r}(\C)$ are in $n$-subgeneral position in 
$\mathbb{P}^{r}(\C)$. For simplicity, we still denote $H_{j}\cap \mathbb{P}^{r}(\C)$ as $H_{j}$ ($1\leq j\leq q$). 
Let $\omega(j)$ be the Nochka weights associated to the hyperplanes $H_{j}$. ($1\leq j\leq q$). 
We denote by $l_{j}$ the corresponding linear form $L_{j}(f)$. For $p\in E$, taking ${\lambda}_{j}=e^{v_{p}(l_{j})}$, 
and using (v) in Theorem \ref{Nochka1}, there exist $L_{p, 1}, \ldots, L_{p, r}$ such that they are linearly independent, 
and that 
\[
\displaystyle \prod_{j=1}^{q} e^{\omega(j){v}_{p}({l}_{j}(f))} \leq \displaystyle \prod_{j=1}^{r} e^{{v}_{p}({l}_{p, j}(f))}\;. 
\] 
This gives
\begin{equation}\label{SMT3-1}
\displaystyle \sum_{j=1}^{q} \omega(j) {v}_{p}({l}_{j}(f))\leq \displaystyle \sum_{j=1}^{r} {v}_{p}({l}_{p, j}(f))\;. 
\end{equation}
For $p\not\in E$, taking ${\lambda}_{j}=e^{v_{p}(l_{j})-\min\{r, v_{p}(l_{j})\}}$, and applying Theorem \ref{Nochka1}, we have 
\begin{equation}\label{SMT3-2}
\displaystyle \sum_{j=1}^{q} \omega(j) [v_{p}(l_{j})-\min\{r, v_{p}(l_{j})\}]\leq \displaystyle \sum_{j=1}^{r} [v_{p}(l_{j})-\min\{r, v_{p}(l_{j})\}]\;.
\end{equation} 
Without loss of generality, we assume that
\begin{equation}\label{SMT3-3}
v_{p}(l_{p, 1})\leq v_{p}(l_{p, 2})\leq \cdots \leq v_{p}(l_{p, r}) \;.
\end{equation}
In the same way of (\ref{Plucker5}), we have
\begin{equation}\label{SMT3-4}
v_{p}(l_{p, j})\leq \delta_{j}(p)\;.
\end{equation}

We consider the case of $p\in E$. In this case, by (\ref{stationary}), (\ref{SMT3-1}) ,and (\ref{SMT3-4}), we have
\[
\displaystyle \sum_{0\leq i\leq r-1}(r-i)\nu_{i}(p) \geq \displaystyle \sum_{j=1}^{r} v_{p}(l_{p, j})-\frac{r(r+1)}{2} 
\geq \displaystyle \sum_{j=1}^{q} \omega(j) v_{p}(l_{j})-\frac{r(r+1)}{2}\;.
\]
Thus 
\begin{equation}\label{SMT3-5}
\displaystyle \sum_{0\leq i\leq r-1} \displaystyle \sum_{p\in E}(r-i)\nu_{i}(p)\geq 
\displaystyle \sum_{j=1}^{q}\sum_{p\in E} \omega(j) v_{p}(l_{j})-\frac{r(r+1)}{2}\sharp E\;.
\end{equation}

We consider the case of $p\not\in E$. Then we can show the following inequality 
\begin{equation}\label{SMT3-6}
\displaystyle \sum_{j=1}^{r}v_{p}(l_{j})-\displaystyle \sum_{j=1}^{r}\min\{r, v_{p}(l_{j})\}\leq \displaystyle \sum_{j=1}^{k} ({\delta}_{i}(p)-i)\;.
\end{equation}
Indeed, assume that  $v_{p}(l_{p, j})\leq r$ for $1\leq j\leq r_{0}$, and $v_{p}(l_{p, j})>r$ for $k_{0}<j\leq r$\:,where 
$1\leq r_{0}\leq r$\:. Then 
\[
\displaystyle \sum_{j=1}^{r}v_{p}(l_{j})-\displaystyle \sum_{j=1}^{r}\min\{r, v_{p}(l_{j})\}=\sum_{j=r_{0}}^{r} (v_{p}(l_{j})-r)\;.
\]
On the other hand, since $\delta_{j}(p)\geq j$ and $v_{p}(l_{p, j})\leq \delta_{j}(p)$ for $j=1,\ldots,k$\;, we get
\[
\displaystyle \sum_{i=0}^{r}(\delta_{i}(p)-i)\geq \displaystyle \sum_{j=r_{0}}^{r}(v_{p}(l_{j})-j)
\geq \displaystyle \sum_{j=r_{0}}^{r} (v_{p}(l_{j})-r)\;.
\]
Combining the above two inequalities, we obtain (\ref{SMT3-6}). 
By (\ref{stationary}) and (\ref{SMT3-6}), we have 
\begin{equation}\label{SMT3-7}
\displaystyle \sum_{0\leq i\leq r-1}\displaystyle \sum_{p\not\in E} (r-i) \nu_{i}(p)= 
\displaystyle \sum_{i=0}^{r} \displaystyle \sum_{p\not\in E} ({\delta}_{i}(p)-i) \geq 
\displaystyle \sum_{j=0}^{r} \displaystyle \sum_{p\not\in E} (v_{p}(l_{j})-\min\{r, v_{p}(l_{j})\})\;.
\end{equation}
By (\ref{SMT3-2}), this implies that 
\begin{equation}\label{SMT3-8}
\displaystyle \sum_{0\leq i\leq r-1}\displaystyle \sum_{p\not\in E} (r-i) \nu_{i}(p)
\geq \displaystyle \sum_{j=0}^{q} \displaystyle \sum_{p\not\in E} \omega(j) [(v_{p}(l_{j})-\min\{r, v_{p}(l_{j})\}]\;.
\end{equation}
By (\ref{stationary}), (\ref{SMT3-5}), and (\ref{SMT3-6}), we get 
\begin{eqnarray*}\label{SMT3-7}
(r+1)\deg(f)+r(r+1)(G-1) =\displaystyle \sum_{p\in \overline{M}} \Bigl(\displaystyle \sum_{0\leq i\leq r-1} (r-i)\nu_{i}(p) \Bigr)  \\
\geq\displaystyle \sum_{j=1}^{q}\omega(j) \deg(f)-\displaystyle \sum_{j=1}^{q} \displaystyle \sum_{p\not\in E} \omega(j)\min\{r, v_{p}(l_{j})\}
-\frac{r(r+1)}{2} \sharp E \\
\end{eqnarray*}
Therefore 
\[
\Bigl(\displaystyle\sum_{j=1}^{q} \omega(j)-(r+1) \Bigr)\deg(f)\leq \displaystyle \sum_{j=1}^{q} \displaystyle \sum_{p\not\in E} \omega(j)\min\{r, v_{p}(l_{j})\}+\frac{r(r+1)}{2}(2G-2+\sharp E)\:.
\]
By (i), (ii) and (iv) in Theorem \ref{Nochka1}, we have
\begin{eqnarray*}
(q-2n+r-1)deg(f)=\theta\Bigl(\displaystyle\sum_{j=1}^{q} \omega(j)-(r+1) \Bigr)\deg(f) \\
\leq \displaystyle \sum_{j=1}^{q} \displaystyle \sum_{p\not\in E} \theta\omega(j)\min\{r, v_{p}(l_{j})\}+\frac{\theta r(r+1)}{2}(2G-2+\sharp E) \\
\leq \displaystyle \sum_{j=1}^{q} \displaystyle \sum_{p\not\in E} \min\{r, v_{p}(l_{j})\}+\frac{r(2n-r+1)}{2}(2G-2+\sharp E)\;\:. \\
\end{eqnarray*}
We have thus proved the theorem.
\end{proof}
\subsection{Ramification estimate and examples}
In the paper \cite{JR}, Jin and Ru obtain the ramification estimate for the Gauss map of algebraic minimal surfaces in 
$\R^{n}$ by using Theorem \ref{SMT3}. We extend the result to pseudo-algebraic minimal surfaces, and get the follwing 
estimate with invariant ``$R$''.
 Here, one says that $g$ is {\rm ramified over a hyperplane} $H=\{[w]\in \mathbb{P}^{n-1}(\C)\:|\: a_{0}w^{0}+ \cdots +a_{n-1}w^{n-1}=0 \}$ {\rm with multiplicity at least $\nu$} 
if all the zeros of the function $g_{H}=(g, A)$ have orders at least $\nu$, where $A=(a_{0},\ldots,a_{n-1})$\;.
If the image of $g$ omits $H$, we shall say that $g$ is {\rm ramified over} $H$ {\rm with multiplicity $\infty$}.
\begin{theorem}[cf.\:Jin and Ru \cite{JR}]\label{S521}
Consider a pseudo-algebraic minimal surface in $\R^{n}$ with the basic domain 
$M=\overline{M}_{G}\backslash \{p_{1},\ldots,p_{k}\}$. 
Let $g\colon M\to \mathbb{P}^{n-1}(\C)$ be its Gauss map and 
$d$ be the degree of $g$ considered as a map $\overline{M}_{G}$. 
Assume that the image $g(M)$ is contained in some $r$-dimensional projective subspace of $\mathbb{P}^{n-1}(\C)$, 
however it is not in any subspace of dimension lower than $r$, where $1\leq r\leq n-1$.
Let $H_{1}, \ldots, H_{q}$ be the hyperplanes in $\mathbb{P}^{n-1}(\C)$, located in general position.
If the map $g$ is ramified over $H_{j}$ with the multiplicity at least ${\nu}_{j}$ for each $j$, then we have
\begin{equation}\label{S522}
\displaystyle \sum_{j=1}^{q} \Bigl(1-\frac{r}{{\nu}_{j}}\Bigr)\leq (2n-r-1)\Bigl(1+\frac{r}{2R}\Bigr),\quad R=\frac{d}{2G-2+k}\geq 1\:.
\end{equation}
In particular, we have
\begin{equation}\label{S523}
\displaystyle \sum_{j=1}^{q} \Bigl(1-\frac{r}{{\nu}_{j}}\Bigr)\leq \frac{(2n-r-1)(r+2)}{2}
\end{equation}
and for algebraic minimal surfaces, the inequality is a strict inequality.
\end{theorem}
\begin{proof}
Consider that the holomorphic $1$-forms ${\phi}_{i}=\partial x^{i}$ ($1\leq i\leq n$). 
For each $j=1,\ldots, k$, let $\mu_{j}$ be the maximum order of poles of ${\phi}_{i}$ at $p_{j}$. 
We can easily find a non-zero vector $(a_{1},\ldots,a_{n})$ so that 
$\phi=a_{1}{\phi}_{1}+\cdots+a_{n}{\phi}_{n}$ has a pole of order $\mu_{j}$ at $p_{j}$.
Applying the Riemann-Roch formula to the meromorphic $1$-form $\phi$ on $\overline{M}_{G}$, we obtain 
\[
d-\displaystyle \sum_{j=1}^{k} {\mu}_{j}=2G-2
\]
Since $\phi$ has poles of order ${\mu}_{j}\geq 1$ by completeness, we get
\begin{equation}\label{S524}
d=2G-2+\sum_{j=1}^{k} {\mu}_{j}\geq 2G-2+k \;,
\end{equation}
and
\begin{equation}
R\geq 1\;.
\end{equation}
When $M$ is an algebraic minimal surfaces, we have ${\mu}_{j}\geq 2$ by the period condition and so $R>1$\;. 

Now, we prove (\ref{S523}). In this situation, the Gauss map $g$ can be extended holomorphically on 
$\overline{M}_{G}$. Let $\{H_{1},\ldots,H_{r_{0}},\hat{H}_{1},\ldots,\hat{H}_{l_{0}}\}$ be the set of totally ramified 
hyperplanes of $g$, located in general position, where $\{H_{1},\ldots,H_{r_{0}}\}$ are exceptional hyperplanes. 
Assume that $g\colon \overline{M}_{G}\to \mathbb{P}^{r}(\C)$ is linearly non-degenerate, where $1\leq r\leq n-1$. 
Apply Theorem\ref{SMT3} to $g$ with $E=\{p_{1},\ldots,p_{k}\}$, we have
\begin{eqnarray*}
(r_{0}+l_{0}-(2n-r-1))d &\leq&  \displaystyle \sum_{j=1}^{r_{0}} \displaystyle \sum_{p \not\in E}\min\{r, v_{p}(L_{j}(g))\} +\displaystyle \sum_{j=1}^{l_{0}} \displaystyle \sum_{p \not\in E}\min\{r, v_{p}(\hat{L}_{j}(g))\} \\
                        &+& \frac{1}{2}r(2n-r-1)\{2(G-1)+k\}\:,
\end{eqnarray*}
where $L_{j}$ are linear forms defining $H_{j}$, $\hat{L}_{j}$ are linear forms defining $\hat{H}_{j}$. 
Since $H_{1},\ldots,H_{r_{0}}$ are exceptional hyperplanes, for $p\not\in E$, $v_{p}(L_{j}(g))=0$ for 
$1\leq j\leq r_{0}$. On the other hand, by the definition, for $p\in M$, we have
\[
\min\{r, v_{p}(\hat{L}_{j}(g))\} \leq r\min\{1, v_{p}(\hat{L}_{j}(g))\} \leq \frac{r}{{\nu}_{j}}v_{p}(\hat{L}_{j}(g))\;.
\] 
Thus we have
\begin{eqnarray*}
(r_{0}+l_{0}-(2n-r-1))d &\leq& \displaystyle \sum_{j=1}^{l_{0}} \displaystyle \sum_{p \not\in E}\frac{r}{{\nu}_{j}}v_{p}(\hat{L}_{j}(g))+\frac{1}{2}r(2n-r-1)\{2(G-1)+k\} \\
                        &=& \displaystyle \sum_{j=1}^{l_{0}} \displaystyle \sum_{p \in M}\frac{r}{{\nu}_{j}}v_{p}(\hat{L}_{j}(g))+\frac{1}{2}r(2n-r-1)\{2(G-1)+k\}\\
                        &=& \displaystyle \sum_{j=1}^{l_{0}} \frac{dr}{\nu_{j}}+\frac{1}{2}r(2n-r-1)\{2(G-1)+k\}\\
\end{eqnarray*}
This implies that
\begin{eqnarray*}
r_{0}+\displaystyle \sum_{j=1}^{l_{0}}\Bigl(1-\frac{r}{{\nu}_{j}}\Bigr) &\leq& (2n-r-1)+\frac{1}{2d}r(2n-r-1)\{2(G-1)+k\} \\
                                                                        &=& (2n-r-1)\Bigl(1+\frac{r}{2R}\Bigr)\;. \\
\end{eqnarray*}
We have thus proved the theorem.
\end{proof}
We have the following result as an immediate consequence of Theorem \ref{S521}.
\begin{corollary}[cf.\:Fujimoto \cite{F3}, Chern-Osserman \cite{CO}, Ru \cite{R1}] \label{coro}
The Gauss map of a non-flat pseudo-algebraic minimal surface in $\R^{n}$ can omit at most 
$n(n+1)/2$ hyperplanes in $\mathbb{P}^{n-1}(\C)$ located in general position. In particular, for an algebraic 
minimal surface, the Gauss map can omit at most $(n-1)(n+2)/2$ hyperplanes in general position in $\mathbb{P}^{n-1}(\C)$.
\end{corollary}

Here, we shall show that, for an arbitrary odd number $n$, the numeber $n(n+1)/2$ of Corollary \ref{coro} is the 
best-possible, namely, there exist pseudo-algebraic minimal surfaces in $\R^{n}$ whose the Gauss maps 
are linearly nondegenerate and omit $n(n+1)/2$ hyperplanes in general position. This result is a modification of Fujimoto's result [\cite{F4}, p193].
\begin{theorem}[cf.\:Fujimoto \cite{F4}]\label{examp}
For an arbitrarily given odd number $n$, there is a pseudo-algebraic minimal surfaces in $\R^{n}$ whose 
Gauss map is linearly nondegenerate and omit $n(n+1)/2$ hyperplanes in $\mathbb{P}^{n-1}(\C)$ located in general position. 
\end{theorem}
To prove this, we first give the following algebraic result.
\begin{lemma}[Fujimoto, \cite{F4} p193]\label{indep}
Let $n$ be an odd number. For $0\leq t\leq (n-1)/2$, we consider $(t+1)n$ polynomials 
\begin{eqnarray*}
f_{i}(u) &=& (u-a_{0})^{n-i} \quad (1\leq i\leq n) \\
f_{n+i}(u) &=& (u-a_{1})^{n-i}(u-b_{1})^{i-1} \quad (1\leq i\leq n) \\
\vdots \\
f_{tn+i}(u) &=& (u-a_{t})^{n-i}(u-b_{t})^{i-1} (1\leq i\leq n)\;, \\
\end{eqnarray*} 
where $a_{\sigma}, b_{\tau}$ are mutually distinct complex numbers. If we take $a_{\sigma}$ and $b_{\tau}$ 
$(0\leq \sigma \leq t$, $1\leq \tau \leq t)$ suitably, then arbitrarily chosen $n$ polynomials among them 
are linearly independent. 
\end{lemma} 
{\em Proof of Theorem \ref{examp}}\quad For a given odd number $n$, we set $m=n-1$, $k=m/2$ and 
define $n$ functions 
\begin{eqnarray*}
h_{2l+1}(z) &=& z^{l}+z^{2k-l} \quad (0\leq l\leq k-1)\;, \\
h_{2l+2}(z) &=& \sqrt{-1}(z^{l}-z^{2k-l}) \quad (0\leq l\leq k-1)
\end{eqnarray*}
and 
\[
h_{2k+1}(z)=2\sqrt{-k}z^{k}\;.
\]
Next, we take suitable constants $a_{\sigma}$ $(0\leq \sigma \leq k)$ and $b_{\tau}$ $(1\leq \tau\leq k)$ such that 
the polynomials $f_{i}$ $(1\leq i\leq q=n(n+1)/2)$ have the properties in Lemma \ref{indep} for $t=k$. By changing 
the variable $u$ suitably, we may assume that $a_{0}=0$. Set
\[
M=\C\backslash \{a_{1},\ldots,a_{k},b_{1},\ldots,b_{k}\}
\] 
and consider the universal covering $\pi\colon \tilde{M}\to M$. Set
\[
\psi(z)=\frac{1}{(z-a_{1})(z-b_{1}) \ldots (z-a_{k})(z-b_{k})}
\]
and define $n$ holomorphic functions $\tilde{g_{i}}=\psi h_{i}$ $(1\leq i\leq n)$\;. Then we see
\[
(\tilde{g_{1}})^{2}+(\tilde{g_{2}})^{2}+\cdots+(\tilde{g_{n}})^{2}=0
\] 
For brevity, we denote the functions $\tilde{g_{i}}\circ\pi$ and $\tilde{g_{i}}$ by the abbreviated notation $g_{i}$ 
in the following.

We consider the functions $x_{i}$ defined by 
\[
x_{i}=\Re \displaystyle \int_{z_{0}}^{z} \phi_{i}
\]
for the holomorphic $1$-forms $\phi_{i}=g_{i}dz$ $(1\leq i\leq n)$. By Enneper-Weierstrass representation 
(Theorem \ref{minimal3}), the surface $x=(x^{1},\ldots,x^{n})\colon \tilde{M}\to \R^{n}$ is a minimal surface. 
Moreover, its Gauss map may be rewritten as $g=(g_{1}:\cdots:g_{n})$ and therefore 
$g=(h_{1}:\cdots:h_{n})$. Since the polynomials
\begin{eqnarray*}
P_{2l+1}(u) &=& u^{l}+u^{2k-l} \quad (0\leq l\leq k-1)\;, \\
P_{2l+2}(u) &=& \sqrt{-1}(u^{l}-u^{2k-l}) \quad (0\leq l\leq k-1) \\
\end{eqnarray*}
and 
\[
P_{2k+1}(u)=2\sqrt{-k}u^{k}
\]
are linearly independent over $\C$, the Gauss map $g$ is linearly nondegenerate. 
Moreover, since $P_{1},\ldots,P_{n}$ give a basis of the vector space of all polynomials of degree $\leq n-1$, 
we can find some constants $c_{ij}$ such that 
\[
f_{i}=\displaystyle \sum_{j=1}^{n} c_{ij}P_{j}\quad (1\leq i\leq q)
\]
Now we consider $q$ hyperplanes 
\[
H_{i}\colon c_{i1}w^{1}+c_{i2}w^{2}+\cdots+c_{in}w^{n}=0 \quad (1\leq i\leq q). 
\]
These are located in general position because $n$ arbitrary polynomials among the $f_{i}'s$ are linearly independent. 
On the other hand, for each $j=1,\ldots,q$ we can write
\[
\displaystyle \sum_{j=1}^{n}c_{ij}h_{j}(z)=\displaystyle \sum_{j=1}^{n} c_{ij}P_{j}=f_{i}(z)=(z-a_{\tau})^{r_{i}}(z-b_{\tau})^{s_{i}}
\]
with suitable non-negative integers $r_{i}, s_{i}$. In view of the definition of $M$, this implies that each $f_{i}(z)$ vanishes 
nowhere on $\tilde{M}$. Consequently, the Gauss map $g$ omits $q$ hyperplanes $H_{j}$ located in 
general position.

The metric on $M$ induced from $\R^{n}$ is given by
\[
ds^{2}=\frac{1}{2} \frac{\sum_{l=0}^{k-1}(|z|^{2l}+|z|^{2(2k-l)}+4k|z|^{2k}}{(|z-a_{1}||z-b_{1}|\cdots|z-a_{k}||z-b_{k}|)^{2}} |dz|^{2}
\] 
and by 
\[
ds^{2}=\frac{1}{2} \frac{\sum_{l=0}^{k-1}(|{\zeta}|^{2l}+|{\zeta}|^{2(2k-l)})+4k|\zeta|^{2k}}{|1-a_{1}\zeta||1-b_{1}\zeta|\cdots |1-a_{k}\zeta||1-b_{k}\zeta|} \frac{|d\zeta|^{2}}{|\zeta|^{2}}
\]
around the $\infty$ if we take a complex local coordinate $\zeta=1/z$. The surface with this metric is complete. 
Indeed, if there is a piecewise smooth curve $\gamma(t)$ $(0\leq t<1)$ in $\tilde{M}$ with finite length, 
which tends to the boundary of $\tilde{M}$, then the curve $\tilde{\gamma}=\pi\circ\gamma$ in $M$ tends to one of the 
points $a_{1},\ldots, a_{k},b_{1},\ldots, b_{k}$ and $\infty$. This is impossible as is easily seen by the above representations of $ds^{2}$. 
We have thus proved the theorem. \quad $\Box$

\end{document}